\documentclass[12pt]{article}
\usepackage{latexsym}
\usepackage{amsfonts}
\usepackage{amsmath}
\usepackage{amsthm}
\usepackage{amssymb}
\usepackage{amscd}
\usepackage[dvips]{graphicx}
\usepackage{color} 
\usepackage{eepic} 
\usepackage{setspace}
\usepackage{calrsfs}

\newdimen\Squaresize \Squaresize=20pt
\newdimen\thickness \thickness=1pt         
                                                    
\def\Square#1{\hbox{\vrule width \thickness
   \vbox to \Squaresize{\hrule height \thickness\vss                            
      \hbox to \Squaresize{\hss#1\hss}
   \vss\hrule height\thickness} 
\unskip\vrule width \thickness} 
\kern-\thickness}                                                            
                               
\def\vsquare#1{\vbox{\Square{$#1$}}\kern-\thickness}
\def\blank{\omit\hskip\Squaresize}

\def\fibyoung#1{\let\\=\cr              
\vbox{\smallskip\offinterlineskip
\halign{&\vsquare{##}\cr #1}}\,}

\def\borderlessrect#1#2{\hbox{\hskip \thickness
   \vbox to \Squaresize{\vskip \thickness \vss
      \hbox to #2 {\hss #1\hss}
   \vss\vskip\thickness} 
\unskip\hskip \thickness} 
\kern-\thickness}                                                            
                               
\def\vborderlessrect#1#2{\vbox{\borderlessrect{$#1$}{#2}}\kern-\thickness}

\def\borderless#1{\omit\vborderlessrect{#1}{\Squaresize}}
\def\borderlessrc#1#2{\omit\vborderlessrect{#1}{#2}}

\def\msquare#1{\vbox{\hbox{\vrule width \thickness
   \vbox to \Squaresize{\hrule height \thickness
      \hbox to \Squaresize{\hfil{\sevenrm #1}}
   \vfil\hrule height\thickness}
\unskip\vrule width \thickness}
\kern-\thickness}\kern-\thickness}

\def\twosquare#1#2{\vbox{\hbox{\vrule width \thickness 
   \vbox to \Squaresize{\hrule height \thickness
      \hbox to \Squaresize{\hfil{\rm #1}}\vss
      \hbox to \Squaresize{\hss{#2}\hss}
   \vfil\hrule height\thickness}
\unskip\vrule width \thickness}
\kern-\thickness}\kern-\thickness}

\def\twoblank#1#2{\vbox{\hbox{
   \vbox to \Squaresize{\vskip 2pt
      \hbox to \Squaresize{\hfil{\sevenrm #1}\ }\vss
      \hbox to \Squaresize{\hss{#2}\hss}
   \vfil}\unskip\kern-\thickness}
}\unskip\kern-\thickness}

\def\young#1{
\def\>{\blank}
\def\<{\borderless}
\def\*{\borderlessrc}
\def\p{\omit\msquare}
\def\t{\omit\twosquare}
\def\b{\omit\twoblank}
\let\\=\cr 
\vbox{\smallskip\offinterlineskip
\halign{&\vsquare{##}\cr #1}}}

\newdimen\smsquaresize \smsquaresize=12pt
\newdimen\smthickness \smthickness=.5pt
\font\smcellfont=cmss8 scaled \magstep0

\def\smsquare#1{\hbox{\vrule width \smthickness
   \unskip\vbox to \smsquaresize{\hrule height \smthickness\vss
      \hbox to \smsquaresize{\hss{\smcellfont #1}\hss}
   \vss\hrule height\smthickness} 
\unskip\vrule width \smthickness} 
\kern-\smthickness}

\def\smvsquare#1{\vbox{\smsquare{$#1$}}\kern-\smthickness}
\def\blank{\omit\hskip\smsquaresize}

\def\smyoung#1{\let\\=\cr 
\vbox{\smallskip\offinterlineskip
\halign{&\smvsquare{##}\cr #1}}}
\newdimen\vsmsquaresize \vsmsquaresize=10pt
\newdimen\vsmthickness \vsmthickness=.5pt
\font\vsmcellfont=cmsl8 scaled \magstep0
\font\vsmletterfont=cmr6 scaled \magstep0

\def\vsmsquare#1{\hbox{\vrule width \vsmthickness
   \unskip\vbox to \vsmsquaresize{\hrule height \vsmthickness\vss
      \hbox to \vsmsquaresize{\hss{\vsmcellfont #1}\hss}
   \vss\hrule height\vsmthickness} 
\unskip\vrule width \vsmthickness} 
\kern-\vsmthickness}
\def\vsmvsquare#1{\vbox{\vsmsquare{#1}}\kern-\vsmthickness}
\def\vsmblank{\omit\hskip\vsmsquaresize}
\def\vsmborderless#1{\hbox{\hskip \vsmthickness\unskip
   \vbox to \vsmsquaresize{\vss
      \hbox to \vsmsquaresize{\hss{\vsmletterfont #1}\hss}
   \vss} 
\unskip\hskip \vsmthickness} 
\kern-\vsmthickness}                                                            \def\vsmvborderless#1{\vbox{\vsmborderless{#1}}\kern-\vsmthickness}

\def\vsmyoung#1{
\def\>{\vsmblank}
\def\<{\omit\vsmvborderless}
\let\\=\cr 
\vbox{\smallskip\offinterlineskip
\halign{&\vsmvsquare{##}\cr #1}}}

\newdimen\edgesize \edgesize=20pt
\newdimen\eedgesize \eedgesize=21pt
\newdimen\doublesize \doublesize=39pt
\newdimen\ddoublesize \ddoublesize=43pt

\newdimen\triplesize \triplesize=62pt
\newdimen\tetrasize \tetrasize=85pt
\newdimen\futosa \futosa=1pt

\def\Hako#1{\hbox{\vrule width \futosa
   \vbox to \edgesize{\hrule height \futosa\vss                            
      \hbox to \edgesize{\hss#1\hss}
   \vss\hrule height\futosa} 
\unskip\vrule width \futosa} 
\kern-\futosa}                                                            
                               
\def\Seihokei#1{\vbox{\Hako{#1}}\kern-\futosa}

\def\Horizontal#1{\hbox{\vrule width \futosa
   \vbox to \edgesize{\hrule height \futosa\vss                            
      \hbox to 41pt{\hss#1\hss}
   \vss\hrule height \futosa} 
\unskip\vrule width \futosa} 
\kern-\futosa}                                                            
                               
\def\Vertical#1{\hbox{\vrule width \futosa
   \vbox to \doublesize{\hrule height \futosa\vss                            
      \hbox to \edgesize{\hss#1\hss}
   \vss\hrule height \futosa} 
\unskip\vrule width \futosa} 
\kern-\futosa}                                                            
                               
\def\NS#1{\vbox{\Hako{$#1$}\vskip20pt}\kern-\futosa}
\def\NSS#1{\vbox{\Hako{$#1$}\vskip41pt}\kern-\futosa}
\def\NSSS#1{\vbox{\Hako{$#1$}\vskip57pt}\kern-\futosa}
\def\H#1{\vbox{\Horizontal{$#1$}}\kern-\futosa}
\def\HH#1#2{\vbox{\Horizontal{$#1$}\Horizontal{$#2$}}\kern-\futosa}
\def\VER#1{\vbox{\Vertical{$#1$}}\kern-\futosa}
\def\VS#1{\vbox{\Vertical{$#1$}\vskip19pt}\kern-\futosa}
\def\VSS#1{\vbox{\Vertical{$#1$}\vskip38pt}\kern-\futosa}
\def\VV#1#2{\vbox{\Vertical{$#1$}\Vertical{$#2$}}\kern-\futosa}
\def\NV#1#2{\vbox{\Hako{$#1$}\Vertical{$#2$}}\kern-\futosa}
\def\NVS#1#2{\vbox{\Hako{$#1$}\Vertical{$#2$}\vskip19pt}\kern-\futosa}
\def\YTT#1#2#3{\vbox{\Horizontal{$#1$}\hbox{\vbox{\Vertical{$#2$}}\kern-\futosa\vbox{\Vertical{$#3$}}\kern-\futosa}}\kern-\futosa}

\def\domino#1{
\def\ns{\omit\NS}
\def\nss{\omit\NSS}
\def\nsss{\omit\NSSS}
\def\h{\omit\H}
\def\hh{\omit\HH}
\def\V{\omit\VER}
\def\Vs{\omit\VS}
\def\Vss{\omit\VSS}
\def\vsss{\omit\VSSS}
\def\vv{\omit\VV}
\def\nv{\omit\NV}
\def\nvs{\omit\NVS}
\def\ytt{\omit\YTT}
\let\\=\cr 
\vbox{\smallskip\offinterlineskip
\halign{&\Seihokei{##}\cr #1}}}

\newdimen\smedgesize \smedgesize=12pt
\newdimen\smeedgesize \smeedgesize=12.5pt
\newdimen\smdoublesize \smdoublesize=24.5pt
\newdimen\smddoublesize \smddoublesize=25pt

\newdimen\smtriplesize \smtriplesize=36.5pt
\newdimen\smtetrasize \smtetrasize=49pt
\newdimen\smfutosa \smfutosa=0.5pt         

\def\smHako#1{\hbox{\vrule width \smfutosa
   \vbox to \smeedgesize{\hrule height \smfutosa\vss                            
      \hbox to \smedgesize{\hss#1\hss}
   \vss\hrule height\smfutosa} 
\unskip\vrule width \smfutosa} 
\kern-\smfutosa}                                                            
                               
\def\smSeihokei#1{\vbox{\smHako{#1}}\kern-\smfutosa}

\def\smHorizontal#1{\hbox{\vrule width \smfutosa
   \vbox to \smeedgesize{\hrule height \smfutosa\vss                            
      \hbox to \smdoublesize{\hss#1\hss}
   \vss\hrule height \smfutosa} 
\unskip\vrule width \smfutosa} 
\kern-\smfutosa}

\def\smVertical#1{\hbox{\vrule width \smfutosa
   \vbox to \smddoublesize{\hrule height \smfutosa\vss
      \hbox to \smedgesize{\hss#1\hss}
   \vss\hrule height \smfutosa} 
\unskip\vrule width \smfutosa} 
\kern-\smfutosa}                                                            
                               
\def\smNS#1{\vbox{\smHako{$#1$}\vskip\smeedgesize}\kern-\smfutosa}
\def\smNSS#1{\vbox{\smHako{$#1$}\vskip\smdoublesize}\kern-\smfutosa}
\def\smNSSS#1{\vbox{\smHako{$#1$}\vskip\smtriplesize}\kern-\smfutosa}
\def\smH#1{\vbox{\smHorizontal{$#1$}}\kern-\smfutosa}
\def\smHH#1#2{\vbox{\smHorizontal{$#1$}\smHorizontal{$#2$}}\kern-\smfutosa}
\def\smHS#1{\vbox{\smHorizontal{$#1$}\vskip\smeedgesize}\kern-\smfutosa}
\def\smVER#1{\vbox{\smVertical{$#1$}}\kern-\smfutosa}
\def\smVS#1{\vbox{\smVertical{$#1$}\vskip\smedgesize}\kern-\smfutosa}
\def\smVSS#1{\vbox{\smVertical{$#1$}\vskip\smdoublesize}\kern-\smfutosa}
\def\smVV#1#2{\vbox{\smVertical{$#1$}\smVertical{$#2$}}\kern-\smfutosa}
\def\smVv#1#2{\vbox{\hbox{\vbox{\smVertical{$#1$}}\kern-\smfutosa\vbox{\smVertical{$#2$}}\kern-\smfutosa}}\kern-\smfutosa}
\def\smNV#1#2{\vbox{\smHako{$#1$}\smVertical{$#2$}}\kern-\smfutosa}
\def\smNVS#1#2{\vbox{\smHako{$#1$}\smVertical{$#2$}\vskip21pt}\kern-\smfutosa}
\def\smYTT#1#2#3{\vbox{\smHorizontal{$#1$}\hbox{\vbox{\smVertical{$#2$}}\kern-\smfutosa\vbox{\Vertical{$#3$}}\kern-\smfutosa}}\kern-\smfutosa}

\def\smdomino#1{
\def\ns{\omit\smNS}
\def\nss{\omit\smNSS}
\def\nsss{\omit\smNSSS}
\def\h{\omit\smH}
\def\hh{\omit\smHH}
\def\hs{\omit\smHS}
\def\V{\omit\smVER}
\def\Vs{\omit\smVS}
\def\Vss{\omit\smVSS}
\def\vsss{\omit\smVSSS}
\def\vv{\omit\smVV}
\def\Vv{\omit\smVv}
\def\nv{\omit\smNV}
\def\nvs{\omit\smNVS}
\def\ytt{\omit\smYTT}
\let\\=\cr 
\vbox{\smallskip\offinterlineskip
\halign{&\smSeihokei{##}\cr #1}}}

\theoremstyle{definition}
\newtheorem{theorem}{Theorem}[section]
\newtheorem{prop}[theorem]{Proposition}
\newtheorem{lemma}[theorem]{Lemma}
\newtheorem{corollary}[theorem]{Corollary}
\newtheorem{definition}[theorem]{Definition}

\newtheorem{conjecture}[theorem]{Conjecture}

\newenvironment{demo}[1]{%
  \trivlist
  \item[\hskip\labelsep
        {\bf #1.}]
}{%
\hfill\qedsymbol
  \endtrivlist
}

\renewcommand{\mathcal}{\mathrsfs}

\newcommand{\p}{\partial}

\renewcommand{\b}{\beta}

\newcommand\Pos{\mathbb{P}}

\renewcommand\tilde{\widetilde}

\newcommand\rdots{\mathinner{\mkern1mu\raise0pt\vbox{\kern7pt\hbox{.}}
     \mkern2mu\raise4pt\hbox{.}\mkern2mu\raise8pt\hbox{.}\mkern1mu}}
\def\TSSCPP#1{\mathcal{T}_{#1}}
\def\TSPP#1{\mathcal{B}_{#1}}
\def\CSPP#1{\mathcal{P}_{#1}}

\def\RDPP#1{\mathcal{D}^\text{R}_{#1}}
\def\CDPP#1{\mathcal{D}^\text{C}_{#1}}
\def\DPP#1{\mathcal{D}_{#1}}

\def\BIJ#1{{\varphi}_{#1}}
\def\GTM#1{{\psi}_{#1}}
\def\DTM#1{{\tau}_{#1}}

\def\BK#1{{\tilde\pi}_{#1}}
\def\FLIP#1{{\pi}_{#1}}
\def\BORDER#1#2#3{{\theta}_{{#1}}{\left({#3}_{\,\geq{#2}}\right)}}

\def\SHAPE#1{{{\operatorname{sh}}\!\left({#1}\right)}}

\def\GCSPP#1{\mathcal{G}_{#1}}

\def\Mod{\operatorname{mod}}

\def\PCSPP#1{\mathcal{Q}_{#1}}
\def\HPCSPP#1{\mathcal{Q}_{#1}^\text{H}}
\def\VPCSPP#1{\mathcal{Q}_{#1}^\text{V}}
\def\diag#1{{\operatorname{diag}\left({#1}\right)}}
\def\Col#1{{\operatorname{Color}\left({#1}\right)}}

\def\V#1#2#3#4{V^{{#1},{#2}}\left({#3};{#4}\right)}

\numberwithin{equation}{section}
\def\rdots{\mathinner{\mkern1mu\raise0pt\vbox{\kern7pt\hbox{.}}
     \mkern2mu\raise4pt\hbox{.}\mkern2mu\raise8pt\hbox{.}\mkern1mu}}

\def\covered{\mathinner{\mkern1mu\raise0pt\vbox{\kern7pt\hbox{$<$}}
     \mkern-4mu\raise2pt\hbox{.}\mkern2mu}}
\def\covers{\mathinner{\mkern1mu\raise0pt\vbox{\kern7pt\hbox{$>$}}
     \mkern-12mu\raise2pt\hbox{.}\mkern8mu}}

\def\defterm#1{{\sl #1}\/}
\def\newterm#1{{\sl #1}\/}

\def\operatorname#1{{\mathrm{#1}\>\!}}

\def\Bbb#1{{\mathbb{#1}}}

\title{
On refined enumerations of
totally symmetric
self-complementary
plane partitions II
}

\author{
Masao Ishikawa\\
\small Faculty of Education, Tottori University\\[-0.8ex]
\small Koyama, Tottori, Japan\\[-0.8ex]
\small \texttt{ishikawa@fed.tottori-u.ac.jp}
}

\date{
\medbreak
}

\begin{document}

\maketitle

\begin{abstract}

In this paper we settle a weak version of a conjecture (i.e. Conjecture~6)
by Mills, Robbins and Rumsey
in the paper ``Self-complementary totally symmetric plane partitions''
{\em J. Combin. Theory Ser. A} {\bf  42}, 277--292.
In other words we show that the number of shifted plane partitions invariant
under the involution $\gamma$ is equal to the number of alternating sign matrices
invariant under the vertical flip.
We also give a determinant expression 
of the general conjecture (Conjecture~6),
but this determinant is still hard to evaluate.
In this paper we introduce two new classes of domino plane partitions,
one has the same cardinality as the set of half-turn symmetric alternating sign matrices
and the other has the same cardinality as the set of vertically symmetric alternating sign matrices.
\end{abstract}

{
\small
\begin{spacing}{0.1}
\tableofcontents 
\end{spacing}
}

\section{
Introduction
}\label{sec:intro}


A totally symmetric self-complementary plane partition is,
by definition,
a plane partition which is invariant under permutation of the three axes 
and which is equal to its complement (cf. \cite{B,I2,MRR2,MRR3,Ro,Sta1}).
The number of totally symmetric self-complementary plane partitions is known to be
the same as that for alternating sign matrices and descending plane partitions. 
But, there are still several interesting conjectures concerning
totally symmetric self-complementary plane partitions (see \cite{MRR2}).
This paper is the succession of my previous paper \cite{I2}
in which we obtain Pfaffian formulae and constant term identities for Conjecture~2, Conjecture~3
and Conjecture~7
as an application of the minor summation formulas of Pfaffians
obtained in \cite{IW1,IW2}.
In this paper we are mainly concerned with  
two other conjectures,
i.e. Conjecture~4 and Conjecture~6,
by Mills, Robbins and Rumsey in the paper \cite{MRR2}.
We will obtain a determinantal formula for Conjecture~6
as an application of essentially a Binet-Cauchy type formula.
We also introduce two new classes of domino plane partitions
which seemingly look closely related to Conjecture~4 and Conjecture~6.

In \cite{MRR2}
Mills, Robbins and Rumsey have introduced a set of 
triangular shifted plane partitions,
which is bijective to the set of totally symmetric self-complementary plane partitions.
In this paper we denote this set by $\TSPP{n}$,
which is defined to be the set of
triangular shifted plane partitions $b=(b_{ij})_{1\leq i\leq j\leq n-1}$
whose parts are $\leq n$,
weakly decreasing along rows and columns,
and all parts in row $i$ are $\geq n-i$.
For $b=(b_{ij})_{1\leq i\leq j\leq n-1}$ in $\TSPP{n}$ and $k\geq1$,
let
\begin{equation}
U_{k}(b)=\sum_{t=1}^{n-k}(b_{t,t+k-1}-b_{t,t+k})
+\sum_{t=n-k+1}^{n-1}\chi\{b_{t,n-1}>n-t\}.
\label{eq:stat_mrr}
\end{equation}
Here we use the convention that $b_{i,n}=n-i$ for all $i$ and $b_{0,j}=n$ for all $j$.
They have also introduced two involutions $\rho$ and $\gamma$ of $\TSPP{n}$ onto itself,
and conjectured that they correspond to the half turn and the vertical flip 
of the alternating sign matrices.
These involutions are defined as follows.
Let $b=(b_{ij})_{1\leq i\leq j\leq n-1}$ be an element of $\TSPP{n}$ 
and let $b_{ij}$ be a part of $b$ off the main diagonal.
Then the \defterm{flip} of the part $b_{ij}$ is the operation of replacing $b_{ij}$ by $b_{ij}'$
where
\begin{equation}
b_{ij}'+b_{ij}=\min(b_{i-1,j},b_{i,j-1})+\max(b_{i,j+1},b_{i+1,j}).
\label{eq:off-diagonal}
\end{equation}
When the part is in the main diagonal,
the flip of a part $b_{ii}$ is the operation replacing $b_{ii}$ by $b_{ii}'$ where
\begin{equation}
b_{ii}'+b_{ii}=b_{i-1,i}+b_{i,i+1}.
\label{eq:main-diagonal}
\end{equation}
An operation $\FLIP{r}$ is defined
to be a map $\TSPP{n}\rightarrow\TSPP{n}$
where $\FLIP{r}(b)$ is the result
of flipping all the $b_{i,i+r-1}$,
$1\leq i\leq n-r$.
We can introduce two involutions $\rho$ and $\gamma$ as follows:
\begin{align}
&\rho=\FLIP{2}\FLIP{4}\FLIP{6}\cdots\\
&\gamma=\FLIP{1}\FLIP{3}\FLIP{5}\cdots
\end{align}
(see (\cite[pp.284,286]{MRR2})).
Mills, Robbins and Rumsey conjectured that the invariants of $\rho$ in $\TSPP{n}$ correspond to the half-turn symmetric
alternating sign matrices,
and the invariants of $\gamma$ in $\TSPP{n}$ correspond to the vertically symmetric
alternating sign matrices.
Let $\TSPP{n}^{\rho}$ (resp. $\TSPP{n}^{\gamma}$) denotes 
the set of elements in $\TSPP{n}$ invariant under $\rho$ (resp. $\gamma$).
These conjectures are stated as follows.
Here,
for the definition of the numbers $A_{n}^\text{HTS}$, $A_{2n+1}^\text{VS}$
and the polynomials $A_{n}^\text{HTS}(t)$, $A_{2n+1}^\text{VS}(t)$,
see the next section.
\begin{conjecture}
\label{conj:HTS}
(\cite[pp.285, Conjecture~4]{MRR2})
Let $n\geq2$ and $r$, $0\leq r< n$, be integers.
Then the number of elements of $\TSPP{n}$ 
with $\rho(b)=b$ and $U_{1}(b)=r$ would be the same as the number of
$n\times n$ alternating sign matrices
invariant under the half turn and satisfying $a_{1,r+1}=1$.
Namely,
 $\sum_{b\in\TSPP{n}^{\rho}}
t^{U_{1}(b)}=A^\text{HTS}_{n}(t)$ would hold.
\end{conjecture}
\begin{conjecture}
\label{conj:VS}
(\cite[pp.286, Conjecture~6]{MRR2})
Let $n\geq1$ be an integer and $r$, $1\leq r\leq 2n-1$, be an integer.
Then the number of elements of $\TSPP{2n+1}$ 
with $\gamma(b)=b$ and $U_{2}(b)=r-1$ would be the same as the number of
$n\times n$ alternating sign matrices with $a_{i1}=1$ and
invariant under the vertical flip.
Namely,
 $\sum_{b\in\TSPP{2n+1}^{\gamma}}
t^{U_{2}(b)}=A^\text{VS}_{2n+1}(t)$ would hold.
\end{conjecture}

In \cite{I2} we have introduced a set $\CSPP{n}$ of column strict plane partitions
which is bijective with the set $\TSPP{n}$ of triangular symmetric plane partitions
(Theorem~\ref{cor:bijection_TSPPtoCSPP}).
Namely $\CSPP{n}$ is,
by definition,
 the set of column-strict plane partitions 
in which each entry in the $j$th column does not exceed $n-j$. 
We recall these plane partitions, the bijections and the statistics in Section~\ref{sec:bijections}.
In Section~\ref{sec:fj},
we translate the involution $\FLIP{r}$  in the words of $\CSPP{n}$,
and find that
the involution $\FLIP{r}$ correspond to a Bender-Knuth type involution $\BK{r}$
which swaps $i$ and $i-1$ in a column-strict plane partition in $\CSPP{n}$.
Let
\begin{align}
&\tilde\rho=\BK{2}\BK{4}\BK{6}\cdots,\\
&\tilde\gamma=\BK{1}\BK{3}\BK{5}\cdots,
\end{align}
and let $\CSPP{n}^{\tilde\rho}$ (resp. $\CSPP{n}^{\tilde\gamma}$) denotes 
the set of invariants of $\tilde\rho$ (resp. $\tilde\gamma$).
Since this Bender-Knuth type involution $\BK{r}$ is, in a sense, ``twisted'' 
(a little different from the ordinary one),
we will see that the set $\CSPP{n}^{\tilde\rho}$ 
is naturally bijective to a set $\GCSPP{n}$ of ``twisted'' domino plane partitions
in Section~\ref{sec:domino}
(see Theorem~\ref{thm:twisted}).
In Section~\ref{sec:gf},
we introduce two new classes of domino plane partitions,
i.e. $\RDPP{n}$ and $\CDPP{n}$.
Namely,
$\RDPP{n}$ (resp. $\CDPP{n}$) is defined to be the set 
of column-strict domino plane partitions whose entries in the $j$th column are $\leq\lceil(n-j)/2\rceil$
and with all rows (resp. columns) of even length.
We also construct a natural bijection 
between $\CSPP{2n+1}^{\tilde\gamma}$ and $\RDPP{2n-1}$
(see Theorem~\ref{thm:bij_domino}).
It seems that $\GCSPP{n}$ and $\CDPP{n}$ ($n\geq1$) have the same number of elements from examples,
but  we don't know how to construct a bijection between them at this point
(see Conjecture~\ref{conj:bijection}).
The following diagram of the bijections give a one-to-one correspondence
between the set $\TSPP{2n+1}^{\gamma}$ of triangular shifted plane partitions invariant under $\gamma$
and the set $\RDPP{2n-1}$ of domino plane partitions:
\[\begin{CD}
\TSPP{2n+1}^{\gamma} @> \BIJ{2n+1} >> \CSPP{2n+1}^{\tilde\gamma} @> \DTM{2n+1} >> \RDPP{2n-1} @> \Phi >> \HPCSPP{2n-1}
\end{CD}\]
Meanwhile, we obtain the following one-to-one correspondence
\[
\TSPP{n}^{\rho} \xrightarrow[]{\BIJ{n}}  \CSPP{n}^{\tilde\rho} \xrightarrow[]{\GTM{n}} \GCSPP{n} \dashrightarrow \CDPP{n} \xrightarrow[]{\Phi} \VPCSPP{n}
\]
between the set $\TSPP{n}^{\rho}$ of triangular shifted plane partitions invariant under $\rho$
and the set $\GCSPP{n}$ of ``twisted'' domino plane partitions,
whereas
we don't know the missing bijection between $\GCSPP{n}$ and $\CDPP{n}$.
In Section~\ref{sec:paired},
we use a plane partition analogue $\Phi$ of the Stanton-White bijection which maps 
a domino plane partition to a paired plane partition
(see \cite{CL,SW}),
that enable us to define a bijection between 
$\RDPP{n}$ (resp. $\CDPP{n}$) of domino plane partitions and
$\HPCSPP{n}$ (resp. $\VPCSPP{n}$) of paired plane partitions.
(see Theorem~\ref{thm:SW}).
Using the generating functions obtained in \cite{I2},
we obtain determinantal formulae for the generating functions
of these sets of paired plane partitions
(see Corollary~\ref{cor:main_gen}).
As a special case we show that Conjecture~\ref{conj:VS}
reduce to the evaluation of the determinant in the following theorem:
\begin{theorem}
\label{thm:VS}
Let  $n\geq2$ be a positive integer.
Let $\det R_{n}^\text{o}(t)=(R_{i,j}^\text{o})_{0\leq i,j\leq n}$
be the $n\times n$ matrix where
\begin{equation*}
R_{i,j}^\text{o}=
\binom{i+j-1}{2i-j}+\left\{\binom{i+j-1}{2i-j-1}+\binom{i+j-1}{2i-j+1}\right\}t+\binom{i+j-1}{2i-j}t^2
\end{equation*}
with the convention that 
$R_{0,0}^\text{o}=R_{0,1}^\text{o}=1$.
Then we obtain
\begin{equation}
\sum_{b\in\TSPP{2n+1}^{\gamma}}
t^{U_{2}(b)}
=\det  R_{n}^\text{o}(t).
\end{equation}
(See Corollary~\ref{cor:result}(ii)).
\end{theorem}
Thus Conjecture~\ref{conj:VS} reduce to prove that $\det D_n(t)=A_{2n+1}^\text{VS}(t)$ would hold.
We also obtain a similar formula for the generating function of $\VPCSPP{n}$
(see Corollary~\ref{cor:result}).
This determinant is also conjectured to be $A^\text{HTS}_{n}(t)$,
but still hard to evaluate
(see Conjecture~\ref{conj:result}).
(About determinant evaluation the reader can consult \cite{Kr2}).
Meanwhile, when $t=1$,
we will find that
we can reduce the evaluation of these determinants to the Andrews-Burge determinant \thetag{\ref{eq:Andrews-Burge}}
(see \cite{A1,A2,AB,MRR3})
and we obtain the result that the number of elements in $\RDPP{n}$ (resp. $\CDPP{n}$) is 
equal to the number $A_{n}^\text{VS}$ (resp. $A_{n}^\text{HTS}$)
of vertically (resp. half-turn) symmetric alternating sign matrices.
Thus we prove Conjecture~\ref{conj:VS} is true when $t=1$,
whereas we can't say it for Conjecture~\ref{conj:HTS} because of the missing bijection.
Anyway we define a new object $\CDPP{n}$ of domino plane partitions
which has the same cardinality with 
the set of half-turn symmetric alternating sign matrices.
In the study of these several classes of plane partitions,
we will see they possess many similarities with Young tableaux
and Schur functions,
but sometimes they are twisted and have mysterious coincidences which we can't explain.

\section{
Preliminaries
}\label{sec:preliminaries}

Let $A_{n}(t)$ be the polynomial defined by
\begin{equation}
A_{n}(t)
=\frac{A_{n}}{\binom{3n-2}{n-1}}
\sum_{r=1}^{n}\binom{n+r-2}{n-1}\binom{2n-1-r}{n-1}t^{r-1},
\end{equation}
where $A_{n}$ is the number defined by
$
A_{n}=\prod_{i=0}^{n-1}\frac{(3i+1)!}{(n+i)!}.
$
It is well-known that $A_{n}$ is the number of alternating sign matrices
and $A_{n}(t)$ is the refined ASM distribution
(see \cite{Ku1,MRR1,Ro,Z2}).
Let $A^\text{HTS}_{n}$ be the number defined by
\begin{equation}
A^\text{HTS}_{2n}=\prod_{i=0}^{n-1}\frac{(3i)!(3i+2)!}{\left\{(n+i)!\right\}^2}
\qquad\text{ and }\qquad
A^\text{HTS}_{2n+1}=
\frac{n!(3n)!}{\left\{(2n)!\right\}^2}\cdot A^\text{HTS}_{2n},
\label{eq:A_HTS}
\end{equation}
which is known to be the number of half-turn symmetric alternating sign matrices (see \cite{Ku3,O2,RS2,Str2}).
The first few terms of \thetag{\ref{eq:A_HTS}} are 
$1$, $2$, $3$, $10$, $25$, $140$, $588$.
We follow \cite{RS2} and define the polynomial $\tilde A^\text{HTS}_{n}(t)$ by
\begin{equation}
\frac{\tilde A^\text{HTS}_{2n}(t)}{\tilde A^\text{HTS}_{2n}}
=\frac{(3n-2)(2n-1)!}{(n-1)!(3n-1)!}
\sum_{r=0}^{n}\frac{\{n(n-1)-nr+r^2\}(n+r-2)!(2n-r-2)!}{r!(n-r)!}t^{r}
\label{eq:tildeAt}
\end{equation}
where $\tilde A^\text{HTS}_{2n}=\prod_{i=0}^{n-1}\frac{(3i)!(3i+2)!}{(3i+1)!(n+i)!}$.
Let
\begin{equation}
A^\text{HTS}_{2n}(t)=\tilde A^\text{HTS}_{2n}(t)A_{n}(t),
\label{eq:HTS_even}
\end{equation}
and
\begin{equation}
A^\text{HTS}_{2n+1}(t)=\frac13\left\{
A_{n+1}(t)\tilde A^\text{HTS}_{2n}(t)+A_{n}(t)\tilde A^\text{HTS}_{2n+2}(t)
\label{eq:HTS_odd}
\right\},
\end{equation}
which is known to be the refined enumeration of half-turn symmetric alternating sign matrices 
weighted by the distribution of one in the top row.
The first few terms of \thetag{\ref{eq:HTS_even}} and \thetag{\ref{eq:HTS_odd}}
are 
$A^\text{HTS}_{2}(t)=1+t$,
$A^\text{HTS}_{3}(t)=1+t+t^2$,
$A^\text{HTS}_{4}(t)=2+3t+3t^2+2t^3$,
$A^\text{HTS}_{5}(t)=3+6t+7t^2+6t^3+3t^4$.

We follow \cite{RS1} and define $A^\text{VS}_{2n+1}$ and $A^\text{VS}_{2n+1}(t)$
as follows.
Let $A^\text{VS}_{2n+1}$ be the number given by
\begin{equation}
A^\text{VS}_{2n+1}
=(-3)^{n^2}\prod_{{1\leq i,j\leq 2n+1}\atop{2|j}}\frac{3(j-i)+1}{j-i+2n+1}
=\frac1{2^n}\prod_{k=1}^{n}\frac{(6k-2)!(2k-1)!}{(4k-2)!(4k-1)!}.
\label{eq:num_VS}
\end{equation}
This number $A^\text{VS}_{2n+1}$ is equal to the number of vertically symmetric alternating sign matrices of size $2n+1$
(see \cite{Ku3,O2,RS1}),
and
the first few terms of \thetag{\ref{eq:num_VS}}
are $1$, $3$, $26$, $646$ and $45885$.
Let  $A^\text{VS}_{2n+1}(t)$ be the polynomial
\begin{equation}
A^\text{VS}_{2n+1}(t)=
\frac{A^\text{VS}_{2n-1}}{(4n-2)!}
\sum_{r=1}^{2n}t^{r-1}
\sum_{k=1}^{r}(-1)^{r+k}
\frac{(2n+k-2)!(4n-k-1)!}{(k-1)!(2n-k)!},
\label{eq:poly_VS}
\end{equation}
which is known to be the refined enumeration of vertically symmetric alternating sign matrices 
weighted by the distribution of one in the first column
(see \cite{RS1}).
For instance,
the first few terms of \thetag{\ref{eq:poly_VS}} are
$A^\text{VS}_{3}(t)=1$,
$A^\text{VS}_{5}(t)=1+t+t^2$,
$A^\text{VS}_{7}(t)=3+6{t}+8{t}^{2}+6{t}^{3}+3{t}^{4}$
and
$A^\text{VS}_{9}(t)=26+78{t}+138{t}^{2}+162{t}^{3}+138{t}^{4}+78{t}^{5}+26{t}^{6}$,
and we have $A^\text{VS}_{2n+1}(1)=A^\text{VS}_{2n+1}$.

\medbreak

Next we recall the terminology of partitions and plane partitions.
We follow the notation in Macdonald \cite{M} and Stanley \cite{Sta2}.
If the reader is familiar with the notion,
he can skip the rest of the section.
Let $\Pos$ denote the set of positive integers.
A partition is a sequence
$\lambda=(\lambda_1,\lambda_2,\dots)$
of non-negative integers in non-increasing order:
$\lambda_1\geq\lambda_2\geq\dots$
and containing only finitely many non-zero terms.
The non-zero $\lambda_i$ are call the \defterm{parts} of $\lambda$.
The number of parts is the \defterm{length} of $\lambda$,
denoted by $\ell(\lambda)$;
and the sum of parts is the \defterm{weight} of $\lambda$,
denoted by $|\lambda|$.
The \defterm{diagram} of a partition $\lambda$ may be formally defined as the set of lattice points $(i,j)\in\Bbb{P}^2$
such that $1\leq j\leq\lambda_i$.
We identify $\lambda$ with its diagram.
The \defterm{conjugate} of a partition $\lambda$ is the partition $\lambda'$
whose diagram is the transpose of the diagram of $\lambda$.
A partition with distinct parts is called a \defterm{strict partition}.
The \defterm{shifted diagram} of a strict partition $\mu$ is the set of lattice points $(i,j)\in\Bbb{P}^2$
such that $i\leq j\leq\mu_i+i$.
We identify a strict partition with its shifted diagram.

A \defterm{plane partition} is an array
$\pi=(\pi_{ij})_{i,j\geq1}$
of nonnegative integers such that $\pi$ has finite support
(i.e. finitely many nonzero entries)
and is weakly decreasing in rows and columns.
If $\sum_{i,j\geq1}\pi_{ij}=n$,
then we write $|\pi|=n$ and say that
$\pi$ is a plane partition of $n$,
or $\pi$ has \defterm{weight} $n$.
A \defterm{part} of a plane partition $\pi=(\pi_{ij})_{i,j\geq1}$
is a positive entry $\pi_{ij}>0$.
The \defterm{shape} of $\pi$ is the ordinary partition $\lambda$
for which $\pi$ has $\lambda_i$ nonzero parts in the $i$th row.
We denote the shape of $\pi$ by $\SHAPE{\pi}$.
We also say that $\pi$ has $r$ \defterm{rows} if $r=\ell(\lambda)$.
Similarly,
$\pi$ has $s$ \defterm{columns} if $s=\ell(\lambda')$.
A plane partition is said to be \defterm{column-strict}
if it is strictly decreasing in columns.

Let $\mu$ be a strict partition.
A \defterm{shifted plane partition} $\tau$ of \defterm{shifted shape} $\mu$ 
is an arbitrary filling of the cells of $\mu$ with nonnegative integers 
such that
each entry is weakly decreasing in rows and columns.
In this paper we allow parts to be zero 
for shifted plane partitions of a fixed shifted shape $\mu$.

\section{
Bijections and Statistics
}\label{sec:bijections}

First we recall the results we obtained in the preceding paper \cite{I2}.
We defined the set $\CSPP{n,m}$ of plane partitions
and studied it intensively.
This set $\CSPP{n,m}$ is also the main object we study in this paper:

\begin{definition}
Let $m$ and $n\geq1$ be nonnegative integers.
Let $\CSPP{n,m}$ denote the set
of column-strict plane partitions $c=(c_{ij})_{1\leq i, j}$
subject to the constraints that
\begin{enumerate}
\item[(C1)]
$c$ has at most $n$ columns;
\item[(C2)]
each part in the $j$th column of $c$ does not exceed $n+m-j$.
\end{enumerate}
An element of $\CSPP{n,m}$ is called a \defterm{restricted column-strict plane partition}.
When $m=0$,
we write $\CSPP{n}$ for $\CSPP{n,0}$.
If a part in the $j$th column of $c$ is equal to $n+m-j$,
we call the part a \defterm{saturated part}.
\end{definition}
Let $c=(c_{ij})_{1\leq i\leq n+m,1\leq j\leq n}$ be a plane partition in $\CSPP{n,m}$
and let $k$ be a positive integer.
Let $c_{\,\geq k}$ denote the plane partition formed by the parts $\geq k$.
Let
\begin{equation}
\BORDER{i}{k}{c}=\sharp\{l:c_{i,l}\geq k\}
\label{eq:border}
\end{equation}
denote the length of the $i$th row of $c_{\,\geq k}$,
i.e.
the rightmost column
 containing a letter $\geq k$ in the $i$th row of $c$.

Let $r$ be an integer such that $1\leq r\leq n+m$.
For $c\in\CSPP{n,m}$
let $\overline U_{r}(c)$ be the number of parts equal to $r$
plus the number of saturated parts less than $r$,
i.e.
\begin{align}
\overline U_{r}(c)=\sharp\{(i,j):c_{ij}=r\}+\sharp\{1\leq k<r:c_{1,n+m-k}=k\}.
\label{eq:statUk}
\end{align}
Especially $\overline U_{1}(c)$ is the number of $1$'s in $c$ 
and $\overline U_{n+m}(c)$ is the number of saturated parts in $c$.
For example,
let
\[
\young{
6&\pmb{\it6}&4&\pmb{\it4}&\pmb{\it3}&1&\pmb{\it1}\\
5&3&3&2&1\\
3&2&2&1\\
1&1\\
}
\]
be an element $c$ of $\CSPP{8}$,
then, the bold faced entries are the saturated parts.
Thus we have 
$\overline U_{1}(c)=6$,
$\overline U_{2}(c)=\overline U_{4}(c)=\overline U_{5}(c)=\overline U_{7}(c)=\overline U_{8}(c)=4$
and
$\overline U_{3}(c)=\overline U_{6}(c)=5$.

\medbreak

We also defined the following set $\TSPP{n,m}$ of shifted
plane partitions in \cite{I2},
which is a generalization of $\TSPP{n}$ defined in \cite[pp.281]{MRR2}.
\begin{definition}
\label{def:TSPP}
(See \cite[Theorem~1]{Kr1}).
Let $m$ and $n\geq1$ be nonnegative integers.
Let $\TSPP{n,m}$ denote the set
of shifted plane partitions $b=(b_{ij})_{1\leq i\leq j}$
subject to the constraints that
\begin{enumerate}
\item[(B1)]
the shifted shape of $b$ is  $(n+m-1,n+m-2,\dots,2,1)$;
\item[(B2)]
$\max\{n-i,0\}\leq b_{ij}\leq n$ for $1\leq i\leq j\leq n+m-1$.
\end{enumerate}
When $m=0$,
we write $\TSPP{n}$ for $\TSPP{n,0}$.
In this paper we call an element of $\TSPP{n,m}$ a \defterm{triangular shifted plane partition}
(abbreviated to TSPP).
\end{definition}
We use the convention that $b_{i,n+m}=n-i$ for all $i$ and $b_{0,j}=n$ for all $j$.
For a $b=(b_{ij})_{1\leq i\leq j\leq n+m-1}$ in $\TSPP{n,m}$ 
and an integer $r$ such that $1\leq r\leq n+m$,
let
\begin{equation}
U_{r}(b)=\sum_{t=1}^{n+m-r}(b_{t,t+r-1}-b_{t,t+r})
+\sum_{t=n+m-r+1}^{n+m-1}\chi\{b_{t,n+m-1}>n-t\}.
\label{eq:stat}
\end{equation}
We put $\overline U_{r}(b)=n+m-1-U_{r}(b)$.
For example,
let
\[
\begin{array}{ccccccc}
8&8&8&8&8&8&8\\
 &8&8&8&8&7&6\\
 & &8&8&7&7&6\\
 & & &6&5&5&4\\
 & & & &4&4&3\\
 & & & & &3&3\\
 & & & & & &1
\end{array}
\]
be an element $b$ of $\TSPP{8}$,
then we have $U_{1}(b)=1$, $U_{2}(b)=U_{4}(b)=U_{5}(b)=U_{7}(b)=U_{8}(b)=3$ and $U_{3}(b)=U_{6}(b)=2$.

In \cite{I2}
we have established a bijection between $\TSPP{n,m}$ and $\CSPP{n,m}$,
and proved that these statistics agree.
\begin{theorem}
\label{cor:bijection_TSPPtoCSPP}
Let $m$ and $n\geq1$ be nonnegative integers
and let $c=(c_{ij})_{1\leq i\leq n+m,1\leq j\leq n}$ be a RCSPP in $\CSPP{n,m}$.
Associate to the array $c=(c_{ij})_{1\leq i\leq n+m,1\leq j\leq n}$ 
the array $b=(b_{ij})_{1\leq i\leq j\leq n+m-1}$ defined by
\begin{equation}
n-b_{ij}=\BORDER{n+m-j}{1-i+j}{c}
\label{eq:CSPPtoTSPP}
\end{equation}
with $1\leq i\leq j\leq n+m-1$.
Then $b$ is in $\TSPP{n,m}$,
 and this mapping $\BIJ{n,m}$,
which associate to a RCSPP $c$ the TSPP $b=\BIJ{n.m}(c)$,
is a bijection of $\CSPP{n,m}$ onto $\TSPP{n,m}$.
Further,
by this bijection,
we have $\overline U_{r}(\BIJ{n,m}(c))=\overline U_{r}(c)$ for any $c\in\CSPP{n,m}$.
\end{theorem}

For instance,
the $b$ and $c$ in the above examples
correspond to each other by this bijection between $\TSPP{8}$ and $\CSPP{8}$.
In \cite[Section~2]{I2}
we defined a set $\TSSCPP{n,m}$ of totally symmetric plane partitions
and constructed the bijections 
$\TSSCPP{n,m}\leftrightarrow\TSPP{n,m}$ 
and $\TSSCPP{n,m}\leftrightarrow\CSPP{n,m}$.
Thus the study of totally symmetric plane partitions
reduce to the study of restricted column-strict plane partitions.

\section{
A twisted Bender-Knuth involution
}\label{sec:fj}

A classical method to prove that a Schur function is symmetric is to define involutions
$s_{i}$ on tableaux which swaps the number of $i$'s and $(i-1)$'s,
for each $i$.
This is well-known as the Bender-Knuth involution (\cite{BK}).
In this section we define a twisted Bender-Knuth involution $\BK{r}$ of the set $\CSPP{n,m}$ of RCSPPs
and show that it correspond to the involution $\FLIP{r}$ of $\TSPP{n,m}$.
We continue to use the convention that $b_{i,n+m}=n-i$ for all $i$ and $b_{0,j}=n$ for all $j$.

In \cite{MRR2}, Mills, Robbins and Rumsey have introduced the notion of flip for $\TSPP{n}$.
We can naturally generalize this notion to $\TSPP{n,m}$ 
by the same equations \thetag{\ref{eq:off-diagonal}} and \thetag{\ref{eq:main-diagonal}},
whereas we have to be careful about the range of $m$.
Let $m$ and $n\geq1$ be non-negative integers.
Let $b=(b_{ij})_{1\leq i\leq j\leq n+m-1}$ be an element of $\TSPP{n,m}$ 
and let $1\leq i<j\leq n+m-1$ so that $b_{ij}$ is a part of $b$ off the main diagonal.
We define the \defterm{flip} of the part $b_{ij}$ as 
the operation of replacing $b_{ij}$ with $b_{ij}'$ with \thetag{\ref{eq:off-diagonal}}.
Note that this operation is always well-defined since $b_{ij}'$ satisfies the axiom \thetag{B2},
and the result of flipping $b_{ij}$ is a shifted plane partition.
In fact $b\in\TSPP{n.m}$ implies that
$\min(b_{i-1,j},b_{i,j-1})\geq b_{ij}\geq\max(b_{i,j+1},b_{i+1,j})$
so that we have $n\geq\min(b_{i-1,j},b_{i,j-1})\geq b_{ij}'\geq\max(b_{i,j+1},b_{i+1,j})\geq0$
and $b_{ij}'\geq b_{i,j+1}\geq n-i$.

When the part is in the main diagonal,
we define the flip of a part $b_{ii}$ as 
the operation replacing $b_{ii}$ with $b_{ii}'$ with \thetag{\ref{eq:main-diagonal}}.
Note that the result of flipping a part $b_{ii}$ in the main diagonal may violate the axiom \thetag{B2} unless $m=0$ or $m=1$.
In fact,
if $m\geq2$,
then $b_{n+m-1,n+m-1}'$ can be negative since $b_{n+m-1,n+m}=1-m<0$.
Thus, hereafter,
we assume $m=0$ or $1$ when we consider a flip of a part in the main diagonal.

Let $1\leq r\leq n+m$ and $b=(b_{ij})_{1\leq i\leq j\leq n+m-1}\in\TSPP{n,m}$.
Define an operation $\FLIP{r}:\TSPP{n,m}\rightarrow\TSPP{n,m}$ by $b\mapsto\FLIP{r}(b)$
where $\FLIP{r}(b)$ is the result of flipping all the $b_{i,i+r-1}$,
$1\leq i\leq n+m-r$.
Since none of these parts of $b$ are neighbors,
the result is independent of the order in which the flips are applied,
and this operation $\FLIP{r}$ is evidently an involution,
i.e. $\FLIP{r}^2=id$.
For instance,
the following TSPP $b\in\TSPP{6}$ is mapped to the following TSPP $\pi_{2}(b)$ by the involution $\pi_2$,
and to the following TSPP $\pi_{1}(b)$ by the involution $\pi_1$;
\[
b=\begin{array}{ccccc}
6&6&6&6&5\\
 &6&5&5&5\\
 & &4&4&4\\
 & & &4&4\\
 & & & &1
\end{array},
\qquad
\pi_{2}(b)=\begin{array}{ccccc}
      6&\pmb{\it6}&      6&      6&      5\\
       &      6&\pmb{\it6}&      5&      5\\
       &       &      4&\pmb{\it4}&      4\\
       &       &       &      4&\pmb{\it2}\\
       &       &       &       &      1
\end{array},
\qquad
\pi_{1}(b)=\begin{array}{ccccc}
\pmb{\it6}&      6&      6&      6&      5\\
       &\pmb{\it5}&      6&      5&      5\\
       &       &\pmb{\it5}&      4&      4\\
       &       &       &\pmb{\it4}&      2\\
       &       &       &       &\pmb{\it4}
\end{array}.
\]

Now we define a Bender-Knuth type involution $\BK{r}:\CSPP{n,m}\rightarrow\CSPP{n,m}$.
This involution $\BK{r}$ is an ``almost Bender-Knuth involution''
which swaps $r$'s and $r-1$'s except the fact that it does not count a saturated $r-1$.
In fact, if it did convert a saturated $r-1$ of $c$ in $\CSPP{n,m}$,
the resulting plane partition could violate the axiom of $\CSPP{n,m}$.
Let see the exact definition.
Let $2\leq r\leq n+m$ and $c\in\CSPP{n,m}$.
Consider the parts of $c$ equal to $r$ or $r-1$.
Since $c$ is column-strict,
some columns of $c$ will contain neither $r$ nor $r-1$,
while some others will contain one $r$ and one $r-1$.
These columns we ignore.
We also ignore an $r-1$ in column $n+m-r+1$,
i.e. we ignore a saturated part which is equal to $r-1$
because a saturated $r-1$ can't be changed to $r$.
The remaining parts equal to $r$ or $r-1$ occur once in each column.
Assume row $i$ has a certain number $k$ of $r$'s followed by a certain number $l$ of $r-1$'s.
Note that we don't count an $r-1$ if it is saturated
so that a saturated $r-1$ always remains untouched.
For example,
the three consecutive rows $i-1$, $i$ and $i+1$ of $c$ could look as follows.
\begin{table}[htb] 
\begin{center}
\begin{tabular}{c||ccc|cccccc|ccc|}
\cline{11-13}
$i-1$& \vdots&        &\multicolumn{2}{c}{\vdots}      &        &      &      &        &      &  $r$   &$\hdots$& $r$  \\\cline{2-4}
$i$  & $r$   &$\hdots$& $r$   & $r$  &$\hdots$&  $r$ & $r-1$&$\hdots$& $r-1$& $r-1$ &$\hdots$& $r-1$\\\cline{11-13}
$i+1$& $r-1$ &$\hdots$& $r-1$ \\
\cline{2-4}
\end{tabular}
\end{center}
\end{table}
In row $i$,
convert the $k$ $r$'s and $l$ $r-1$'s to $l$ $r$'s and $k$ $r-1$'s.
It is easy to see that the resulting array satisfies the axioms \thetag{C1} and \thetag{C2} of $\CSPP{n,m}$.
Define an operation $\BK{r}:\CSPP{n,m}\rightarrow\CSPP{n,m}$ by $c\mapsto\BK{r}(c)$
where $\BK{r}(c)$ is the result of swapping $r$'s and $r-1$'s in row $i$ of $c$ by this twisted rule for $1\leq i\leq n+m-r$.
For example,
if $n=6$, $m=0$ and $r=2$,
then the left below RCSPP $c$ corresponds to the right below RCSPP $\BK{2}(c)$
by $\BK{2}$.
\[
\hbox to 40pt{\vbox{$c=$\vskip40pt}}
\young{
\pmb{\it5}&3&1&1&\pmb{\it1}\\
3&2\\
2&1\\
}
\qquad
\hbox to 60pt{\vbox{$\BK{2}(c)=$\vskip40pt}}
\young{
\pmb{\it5}&3&2&\pmb{\it2}&\pmb{\it1}\\
3&2\\
1&1\\
}
\]

Next assume $m=0$ or $1$, and let $c$ be a plane partition in $\CSPP{n,m}$.
Set $\lambda_i=\BORDER{i}{2}{c}$ to be the number of parts $\geq2$ in the $i$th row of $c$.
Assume the $i$th row contains a certain number $k$ of $1$'s followed by a certain number $l$ of blank positions
which 1's can be put in,
so that we have $k+l=n+m-1-\lambda_1$ if $i=1$, $k+l=\lambda_{i-1}-\lambda_i$ otherwise.
Change the number of $1$'s from $k$ to $l$ in row $i$ for $1\leq i\leq n+m-1$.
It is also easy to see that the resulting array,
say $\BK{1}(c)$,
satisfies the constraints \thetag{C1} and \thetag{C2}.
For example,
if $c$ is as above,
then $\BK{1}(c)$ is as follows:
\[
\hbox to 60pt{\vbox{$\BK{1}(c)=$\vskip40pt}}
\young{
\pmb{\it5}&3\\
3&2\\
2\\
1\\
}.
\]
This mapping $\BK{1}:\CSPP{n,m}\rightarrow\CSPP{n,m}$ is well-defined for $m=0,1$,
and is evidently an involution.
We call this involution $\BK{r}$,
$1\leq r\leq n+m$,
a \defterm{twisted Bender-Knuth involution} (abbreviated to the \defterm{TBK involution}).
Note that when $r=n+m$, $\FLIP{n+m}$ and $\BK{n+m}$ are both the identity mapping since there are no parts affected
by the operations.
The following proposition
corresponds to Theorem~2 of \cite[pp.283]{MRR2},
whereas there is no need of proof 
since it is clear from the above definition.
\begin{prop}
\label{prop:distribution}
Let $m$ and $n\geq1$ be non-negative integers.
Let $2\leq r\leq n+m$
and let $c$ in $\CSPP{n,m}$.
Then
\begin{equation*}
\overline U_{r}\left(\BK{r}(c)\right)=\overline U_{r-1}\left(c\right)
\text{ and }
\overline U_{r}\left(c\right)=\overline U_{r-1}\left(\BK{r}(c)\right)
\ \Box
\end{equation*}
\end{prop}
The following theorem tells us that the involution $\FLIP{r}$ of $\TSPP{n,m}$ corresponds to $\BK{r}$ of $\CSPP{n.m}$
if we identify $\TSPP{n,m}$ with $\CSPP{n.m}$ by the bijection $\BIJ{n,m}$ 
defined in Corollary~\ref{cor:bijection_TSPPtoCSPP}.
\begin{theorem}
\label{th:flip}
Let $m$ and $n\geq1$ be non-negative integers and let $1\leq r\leq n+m$.
Assume $m=0$ or $1$ if $r=1$.
Then we have
\begin{equation*}
\FLIP{r}\left(\BIJ{n,m}\left(c\right)\right)
=\BIJ{n,m}\left(\BK{r}(c)\right).
\end{equation*}
\end{theorem}

In \cite[pp.284]{MRR2},
Mills, Robbins and Rumsey defined an involution
$\rho$ of $\TSPP{n}$ by
$
\rho=\FLIP{2}\FLIP{4}\cdots
$
and presented a conjecture 
(Conjecture~\ref{conj:HTS})
that this involution $\rho$ corresponds 
to the half turn of an alternating matrix.
We naturally generalize this definition to $\TSPP{n,m}$
and use the same symbol for
the involution $\rho:\TSPP{n,m}\rightarrow\TSPP{n,m}$ defined by
\begin{equation}
\rho=\pi_{2}\pi_{4}\cdots
\end{equation}
where the product is over all $\pi_{i}$ with $i$ even and $\leq n$.
Let $\TSPP{n,m}^{\rho}$ denote
the set of elements of $\TSPP{n,m}$ invariant under $\rho$,
i.e. $\TSPP{n,m}^{\rho}=\{b\in\TSPP{n,m}\,|\,\rho(b)=b\}$.
For instance,
the following TSPP is an element of $\TSPP{8}^{\rho}$:
\[
\begin{array}{ccccccc}
      8&\pmb{\it8}&      8&\pmb{\it8}&      7&\pmb{\it7}&      7\\
       &      8&\pmb{\it8}&      8&\pmb{\it7}&      7&\pmb{\it7}\\
       &       &      7&\pmb{\it7}&      7&\pmb{\it7}&      7\\
       &       &       &      6&\pmb{\it6}&      6&\pmb{\it5}\\
       &       &       &       &      6&\pmb{\it5}&      4\\
       &       &       &       &       &      4&\pmb{\it3}\\
       &       &       &       &       &       &      1
\end{array}
\]
By Theorem~\ref{th:flip},
we can reduce the properties of $\rho$ 
to those of the corresponding involution
$\tilde\rho:\CSPP{n,m}\rightarrow\CSPP{n,m}$ defined by
\begin{equation}
\tilde\rho=\BK{2}\BK{4}\cdots
\end{equation}
where the product is over all $\BK{i}$ with $i$ even and $\leq n$.
In other words,
given a plane partition $c$ in $\CSPP{n,m}$,
$\tilde\rho$ swaps $1$'s and $2$'s in $c$ by the TBK involution,
then swap $3$'s and $4$'s in $c$ 
and so on.
The resulting plane partition does not depend on the order of the swaps,
and is an element of $\CSPP{n,m}$.
Let $\CSPP{n,m}^{\tilde\rho}$ denote
the set of elements of $\CSPP{n,m}$ which is invariant under $\tilde\rho$,
i.e. $\CSPP{n,m}^{\tilde\rho}=\{c\in\CSPP{n,m}\,|\,\tilde\rho(c)=c\}$.
For example,
if $n=8$,
then the following RCSPP in $\CSPP{8}$ is invariant under $\tilde\rho$:
\[
\young{
\pmb{\it7}&4&4&3&2&1&\pmb{\it1}\cr
6&3&2&1\cr
5&2\cr
2&1\cr
1\cr
}
\]
Thus
$\CSPP{1}^{\tilde\rho}=\{\emptyset\}$,
$\CSPP{2}^{\tilde\rho}=\left\{\emptyset,\,\vbox to 10pt{\vss\smyoung{\pmb{\it1}\\}\vss}\,\right\}$,
and $\CSPP{3}^{\tilde\rho}$ is composed of the following 3 RCSPPs:
\begin{equation*}
\hbox to 20pt{\vbox{$\emptyset$\vskip12pt}}
\qquad\qquad
\smyoung{
\pmb{\it2}\\
1\\
}
\qquad\qquad
\smyoung{
\pmb{\it2}&\pmb{\it1}\\
1\\
}
\end{equation*}
$\CSPP{4}^{\tilde\rho}$ 
is composed of the following 10 elements:
\begin{equation*}
\hbox to 30pt{\vbox{$\emptyset$\vskip15pt}}
\qquad\qquad
\vbox{\smyoung{2&1\\}\vskip12pt}
\qquad\qquad
\vbox{\smyoung{2&1&\pmb{\it1}\\}\vskip12pt}
\qquad\qquad
\smyoung{2\\1\\}
\qquad\qquad
\smyoung{2&\pmb{\it2}\\1&1\\}
\end{equation*}
\begin{equation*}
\vbox{\smyoung{2&\pmb{\it2}&\pmb{\it1}\\1&1\\}\vskip11pt}
\qquad\qquad
\vbox{\smyoung{\pmb{\it3}\\}\vskip22pt}
\qquad\qquad
\smyoung{\pmb{\it3}\\2\\1\\}
\qquad\qquad
\smyoung{\pmb{\it3}&\pmb{\it2}\\2&1\\1\\}
\qquad\qquad
\smyoung{\pmb{\it3}&\pmb{\it2}&\pmb{\it1}\\2&1\\1\\}
\end{equation*}
$\CSPP{5}^{\tilde\rho}$ has 25 elements,
and $\CSPP{6}^{\tilde\rho}$ has 140 elements.

Also in \cite[pp.286]{MRR2}
the involution $\gamma=\FLIP{1}\FLIP{3}\FLIP{5}\cdots$ on $\TSPP{n}$
is defined and
conjectured to have the same effect as the flip of an alternating matrix
around the vertical axis
(Conjecture~\ref{conj:VS}).
Naturally we can generalize this definition to the involution $\gamma:\TSPP{n,m}\rightarrow\TSPP{n,m}$, $m=0,1$,
defined by
\begin{equation}
\gamma=\pi_{1}\pi_{3}\pi_{5}\cdots
\end{equation}
where the product is over all $\pi_{i}$ with $i$ odd and $\leq n$.
Let $\TSPP{n,m}^{\gamma}$ denote the set of elements of $\TSPP{n,m}$
invariant under $\gamma$.
For instance,
the following TSPP in $\TSPP{7}$ is invariant under $\gamma$:
\[
\begin{array}{cccccc}
\pmb{\it7}&      7&\pmb{\it7}&      7&\pmb{\it7}&      7\\
       &\pmb{\it6}&      5&\pmb{\it5}&      5&\pmb{\it5}\\
       &       &\pmb{\it5}&      5&\pmb{\it5}&      5\\
       &       &       &\pmb{\it5}&      5&\pmb{\it4}\\
       &       &       &       &\pmb{\it4}&      3\\
       &       &       &       &       &\pmb{\it2}
\end{array}.
\]
The corresponding involution $\tilde\gamma:\CSPP{n,m}\rightarrow\CSPP{n,m}$, $m=0,1$,
is defined by
\begin{equation}
\tilde\gamma=\BK{1}\BK{3}\BK{5}\cdots
\end{equation}
where the product is over all $\BK{i}$ with $i$ odd and $\leq n$.
For $m=0,1$,
let $\CSPP{n,m}^{\tilde\gamma}$ denote the set of elements of $\CSPP{n,m}$
invariant under $\tilde\gamma$,
i.e. $\CSPP{n,m}^{\tilde\gamma}=\{c\in\CSPP{n,m}\,|\,\tilde\gamma(c)=c\}$.
But,
since $\CSPP{n,1}=\CSPP{n+1,0}$ which implies $\CSPP{n,1}^{\tilde\gamma}=\CSPP{n+1}^{\tilde\gamma}$,
we only need to study $\CSPP{n}^{\tilde\gamma}$.
Also note that $\CSPP{n}^{\tilde\gamma}=\emptyset$ unless $n$ is odd.
In fact, for $c\in\CSPP{n,m}$, there are, in total,
 exactly $n-1$ positions where one can put $1$'s.
But, if $c$ is invariant under $\tilde\gamma$,
then those positions must be half filled.
Further one can easily see that the shape of $c_{\geq2}$ must be even.
For example,
if $n=7$,
the following RCSPP in $\CSPP{7}$ is invariant under $\tilde\gamma$.
\[
\young{
5&5&3&2&1\cr
4&4&1\cr
3&3\cr
2&2\cr
1\cr
}
\]
Thus we have 
$\CSPP{3}^{\tilde\gamma}=\left\{\,\vbox{\smyoung{1\\}}\,\right\}$,
$\CSPP{5}^{\tilde\gamma}$ is composed of the following 3 RCSPPs:
\begin{equation*}
\vbox{\smyoung{1&1\\}\vskip22pt}
\qquad\qquad
\vbox{\smyoung{3&2&1\\1\\}\vskip11pt}
\qquad\qquad
\vbox{\smyoung{3&3&1\\2&2\\1\\}}
\end{equation*}
and $\CSPP{5}^{\tilde\gamma}$ has 26 elements.

By Theorem~\ref{th:flip},
we have established a bijection 
between $\TSPP{n,m}^{\rho}$ and $\CSPP{n,m}^{\tilde\rho}$
and a bijection 
between $\TSPP{n,m}^{\gamma}$ and $\CSPP{n,m}^{\tilde\gamma}$.
By this bijection the weight function $\overline U_{r}$ on $\TSPP{n,m}$
is exactly the same as $\overline U_{r}$ on $\CSPP{n,m}$.
Thus we study $\CSPP{n,m}^{\tilde\rho}$ and $\CSPP{n,m}^{\tilde\gamma}$
with this weight $\overline U_{r}$ in the rest of this paper.

\begin{demo}{Proof of Theorem~\ref{th:flip}}
Let $c=(c_{ij})_{1\leq i\leq n+m-1,1\leq j\leq n}$ be a plane partition in $\CSPP{n,m}$.
We set $b=(b_{ij})_{1\leq i\leq j\leq n+m-1}=\BIJ{n,m}(c)\in\TSPP{n,m}$
to be the TSPP mapped by the bijection
and set $b'=(b_{ij}')_{1\leq i\leq j\leq n+m-1}=\FLIP{r}(b)\in\TSPP{n,m}$
to be the flipped result of $b$.
Let $c'=(c_{ij}')_{1\leq i\leq n+m-1,1\leq j\leq n}=\BIJ{n,m}^{-1}(b')\in\CSPP{n,m}$ be 
the corresponding RCSPP.

First, assume $2\leq r\leq n+m$.
Then $b'$ is the resulting TSPP
of flipping the part $b_{i,i+r-1}$,
$i=1,\dots,n+m-r$,
i.e.
replacing $b_{i,i+r-1}$ by $b_{i,i+r-1}'$ for $i=1,\dots,n+m-r$
where
\begin{equation*}
b_{i,i+r-1}'+b_{i,i+r-1}
=\min(b_{i-1,i+r-1},b_{i,i+r-2})+\max(b_{i,i+r},b_{i+1,i+r-1}).
\end{equation*}
If we restate this operation $\FLIP{r}:b\mapsto b'$ by the bijection rule \thetag{\ref{eq:CSPPtoTSPP}},
then this corresponds to the following operation $\BK{r}:c\mapsto c'$:
\begin{enumerate}
\item[(i)]
$\BORDER{i}{r}{c'}$ is given by 
\begin{align*}
\BORDER{1}{r}{c'}+\BORDER{1}{r}{c}
=\max(\BORDER{1}{r+1}{c},\BORDER{2}{r-1}{c})
+\min(n+m-r,\BORDER{1}{r-1}{c})
\end{align*}
if $i=1$,
and
\begin{align*}
\BORDER{i}{r}{c'}+\BORDER{i}{r}{c}
=\max(\BORDER{i}{r+1}{c},\BORDER{i+1}{r-1}{c})
+\min(\BORDER{i-1}{r+1}{c},\BORDER{i}{r-1}{c})
\end{align*}
if $i=2,\dots,n+m-r$.
\item[(ii)]
If $k\neq r$,
then
$\BORDER{i}{k}{c'}=\BORDER{i}{k}{c}$ for all $i$.
\end{enumerate}
This means the operation $\BK{r}$ changes the number of the letters $r$'s and $r-1$'s 
in each row $i$ of $c$,
and keeps other letters invariant.
If one inspects these rules carefully,
then these rules gives precisely the TBK involution.

Next, assume $r=1$ and $m=0,1$.
Then $b'$ is the resulting TSPP of 
replacing $b_{ii}$ by $b_{ii}'$ for $i=1,\dots,n+m-1$
where
\begin{equation*}
b_{ii}'+b_{ii}
=b_{i-1,i}+b_{i,i+1}.
\end{equation*}
Using the rule \thetag{\ref{eq:CSPPtoTSPP}} again,
the operation $\FLIP{r}:b\mapsto b'$ 
corresponds to the operation $\BK{r}:c\mapsto c'$ with
$\BORDER{i}{r}{c'}$ given by 
\begin{align*}
\BORDER{1}{1}{c'}+\BORDER{1}{1}{c}
=\BORDER{1}{2}{c}+n+m-1
\end{align*}
if $i=1$,
and
\begin{align*}
\BORDER{i}{1}{c'}+\BORDER{i}{1}{c}
=\BORDER{i}{2}{c}+\BORDER{i-1}{2}{c}
\end{align*}
if $i=2,\dots,n+m-1$.
This is exactly the operation $\BK{1}:c\mapsto c'$ we defined above.
This completes the proof.
\end{demo}

\section{
Twisted domino plane partitions
}\label{sec:domino}

In this section
we consider the invariants
of the involution $\tilde\rho$ defined in the previous section,
which leads us to define a notion of generalized domino plane partitions.
In fact,
here,
we consider a natural bijection as follows.
Assume $c\in\CSPP{n,m}$ is invariant under $\tilde\rho$,
i.e. invariant under $\BK{2i}$ for any $i$.
Then we replace paired $2i$ and $2i-1$ in a column by a domino,
and replace non-paired $k$ $2i$'s and $k$ $2i-1$'s in a row by $k$ dominoes.
But there may remain some saturated parts unchanged,
thus we obtain a ``twisted'' domino plane partition by this operation.
In this manner we give a natural bijection of $\CSPP{n,m}^{\tilde\rho}$ to a new object
which we denote by $\GCSPP{n,m}$.
But we don't know how to count this set at this point.
Anyway let's start with definitions.

\medbreak

A \defterm{domino} is a special kind of skew shape consists of two squares.
A $1\times2$ domino is called a \defterm{horizontal domino}
while a $2\times1$ domino is called a \defterm{vertical domino}.
Let $\lambda$ be a partition.
A \defterm{generalized domino plane partition of shape $\lambda$}
consists of a tiling of the shape $\lambda$ by means of 
ordinary $1\times1$ squares and dominoes,
and a filling of each square or domino with
a positive integer so that
the integers are weakly decreasing along either rows or columns.
The integers in the squares or dominoes are called \defterm{parts}.
In this paper we call a part a \defterm{single} part
if it is in a square,
i.e. not in a domino.
Further we call it
a \defterm{domino plane partition}
if the shape $\lambda$ is tiled with only dominoes,
i.e. without a single square.
We say that a part is in the $i$th row (resp. $j$th column)
if the square or domino which contains the number
intersects with the $i$th row (resp. $j$th column) of $\lambda$.
A (generalized) domino plane partition is said to be \defterm{column-strict}
if it is strictly decreasing along each column.
For example,
the left-below is a column-strict generalized domino plane partition of shape $(4,3,2,1)$,
while
the right-below is a column-strict domino plane partition of shape $(4,4,2)$.
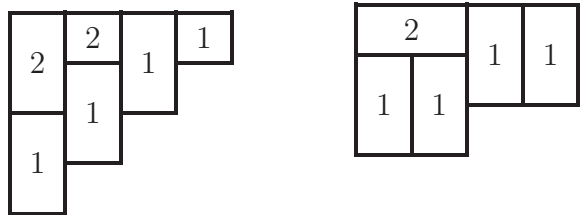
\[
\domino{
\vv{2}{1}&\nvs{2}{1}&\Vss{1}&\nsss{1}\cr
}
\qquad\qquad
\vbox{\domino{
\ytt{2}{1}{1}&\Vs{1}&\Vs{1}
\cr
}\vskip22pt}
\]
\begin{definition}
Let $m$ and $n\geq1$ be nonnegative integers.
Let $\GCSPP{n,m}$
denote the set of column-strict generalized domino plane partitions $c$ subject to the constraints that
\begin{enumerate}
\item[(E1)]
$c$ has at most $n$ columns;
\item[(E2)]
each part in the $j$th column does not exceed $\lceil(n+m-j)/2\rceil$;
\item[(E3)]
If a part in the $j$th column is equal to $\lceil(n+m-j)/2\rceil$
and $n+m-j$ is odd,
then it must be a single part.
On the other hand,
if a single part appears in the $j$th column of $c$,
then $n+m-j$ must be odd and it must be equal to $\lceil(n+m-j)/2\rceil$.
\end{enumerate}
We call an element in $\GCSPP{n,m}$
a \defterm{twisted domino plane partition},
and we simply write $\GCSPP{n}$ for $\GCSPP{n,0}$.
If $c$ is in $\GCSPP{n,m}$,
we call a single part in $c$ a \defterm{saturated} part,
which can appear only in the first row.
Thus a saturated part is always a single part and vice versa,
which equals $\lceil(n+m-j)/2\rceil$ appearing in the cell $(1,j)$
where $n+m-j$ is odd.
\end{definition}
For instance,
$\GCSPP{1}=\{\emptyset\}$,
$\GCSPP{2}=\left\{\emptyset,\,\vbox{\vss\smyoung{\pmb{\it1}\\}\vss}\,\right\}$,
$\GCSPP{3}$
is composed of the following 3 elements:
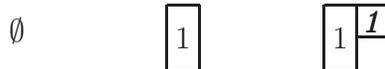
\[
\hbox to 30pt{\vbox{$\emptyset$\vskip10pt}}
\qquad\qquad
\smdomino{\V{1}\cr}
\qquad\qquad
\smdomino{\V{1}&\ns{\pmb{\it1}}\cr}
\]
$\GCSPP{4}$
is composed of the following 10 elements:
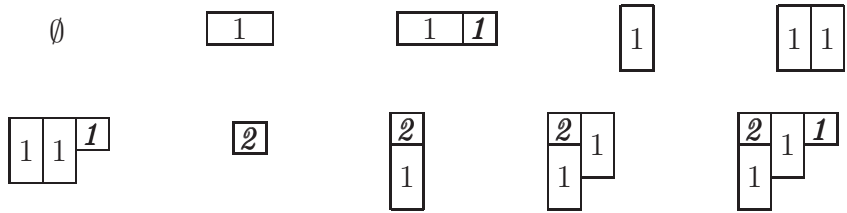
\begin{equation*}
\hbox to 30pt{\vbox{$\emptyset$\vskip10pt}}
\qquad\qquad
\vbox{\smdomino{\h{1}\cr}\vskip10pt}
\qquad\qquad
\vbox{\smdomino{\h{1}&$\pmb{\it1}$\cr}\vskip10pt}
\qquad\qquad
\smdomino{\V{1}\cr}
\qquad\qquad
\smdomino{\V{1}&\V{1}\cr}
\end{equation*}
\begin{equation*}
\vbox{\smdomino{\V{1}&\V{1}&\ns{\pmb{\it1}}\cr}\vskip10pt}
\qquad\qquad
\vbox{\smdomino{$\pmb{\it2}$\cr}\vskip21pt}
\qquad\qquad
\smdomino{$\pmb{\it2}$\cr\V{1}\cr}
\qquad\qquad
\smdomino{\nv{\pmb{\it2}}{1}&\Vs{1}\cr}
\qquad\qquad
\smdomino{\nv{\pmb{\it2}}{1}&\Vs{1}&\nss{\pmb{\it1}}\cr}
\end{equation*}
$\GCSPP{5}$ has 25 elements
and $\GCSPP{6}$ has 140 elements.
For example,
if $c\in\GCSPP{5}$,
 all parts in the 1st and 2nd columns are $\leq2$,
and all parts in the 3rd and 4th columns are $\leq1$.
A saturated (=single) part equal to $2$ can appear in the cell $(1,2)$,
and a saturated part equal to $1$ can appear in the cell $(1,4)$.

For $c\in\GCSPP{n,m}$,
let $\overline U_{1}(c)$ denote the number of $1$'s in $c$.
From the above examples
we obtain that the first few terms of $\sum_{c\in\GCSPP{r}}t^{\overline U_{1}(c)}$
are $1$, $1+t$, $1+t+t^2$,
$(1+t)(2+t+2t^2)$, $3+6t+7t^2+6t^3+3t^4$ and $5(1+t)(1+t^2)(2+3t+2t^2)$.

\bigbreak

Let $m$ and $n\geq1$ be nonnegative integers.
Let's consider the situation when we apply $\tilde\pi_{r}$
to a plane partition $c$ in $\CSPP{n,m}$.
We say an entry $r$ or $r-1$ in $c$ is \newterm{free}
if it is not saturated $r-1$ nor
there is no corresponding $r-1$ or $r$ in the same column.
If $r$ and $r-1$ are in the same column,
we say they are \newterm{paired}.
Thus an $r-1$ in $c$ can be free, paired or saturated,
whereas an $r$ in $c$ can be free or paired.

Let $c$ be a plane partition in $\CSPP{n,m}^{\tilde\rho}$,
i.e. $c$ is invariant under each $\tilde\pi_{2r}$ for $r\geq1$.
We associate to $c$
a column-strict generalized domino plane partition $d$ as follows.
For each $r\geq1$, we replace paired $2r$ and $2r-1$ by a vertical domino containing $r$,
and, if row $i$ contains $k$ free $2r$'s and $k$ free $2r-1$'s,
then we replace these with $k$ horizontal dominoes containing $r$.
Finally a saturated part equal to $2r-1$ should be replaced by a single part $r$.
Let us denote by $\GTM{n,m}(c)$ the resulting generalized domino plane partition $d$.
We may use the abbreviated notation $\GTM{n}(c)$ for $\GTM{n,0}(c)$.
For example,
the left-below plane partition in $\CSPP{8}^{\tilde\rho}$
is mapped to the right-below generalized plane partition in $\GCSPP{8}$
by $\GTM{8}$:
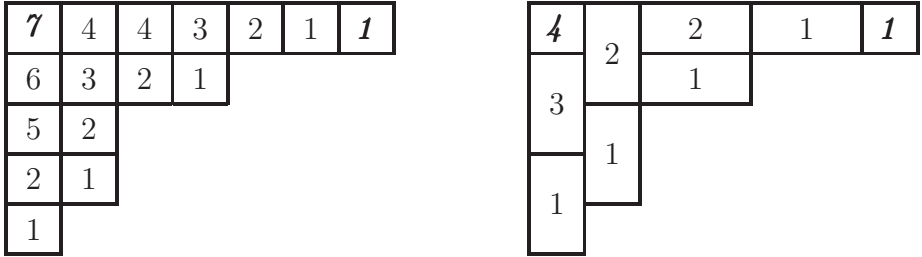
\begin{equation*}
\young{
\pmb{\it7}&4&4&3&2&1&\pmb{\it1}\cr
6&3&2&1\cr
5&2\cr
2&1\cr
1\cr
}
\qquad\qquad
\hbox{
\vbox {
\Hako{$\pmb{\it4}$}
\kern-\futosa
\Vertical{$3$}
\kern-\futosa
\Vertical{$1$}
}\kern-\futosa
\vbox{
\Vertical{$2$}
\kern-\futosa
\Vertical{$1$}
\vskip19pt
}\kern-\futosa
\vbox{
\Horizontal{$2$}
\kern-\futosa
\Horizontal{$1$}
\vskip57pt
}\kern-\futosa
\vbox{
\Horizontal{$1$}
\vskip76pt
}\kern-\futosa
\vbox{
\Hako{$\pmb{\it1}$}
\vskip76pt
}
}
\end{equation*}

%
%
%
\begin{theorem}
\label{thm:twisted}
Let $m$ and $n\geq1$ be nonnegative integers
and let $c=(c_{ij})_{1\leq i\leq n+m,1\leq j\leq n}$ be a plane partition
in $\CSPP{n,m}^{\tilde\rho}$.
Associate to $c$ the generalized domino plane partition $\GTM{n,m}(c)$
as above.
Then $\GTM{n,m}(c)$ is in $\GCSPP{n,m}$,
and $\GTM{n,m}$ is a bijection between $\CSPP{n,m}^{\tilde\rho}$ and $\GCSPP{n,m}$.
\end{theorem}

\section{
Domino plane partitions
}\label{sec:gf}

In this section,
we introduce two important classes of domino plane partitions,
i.e. $\RDPP{n,m}$ and $\CDPP{n,m}$.
The main result of this section is Theorem~\ref{thm:bij_domino}
which shows that there is a bijection between 
the set $\CSPP{2n+1}^{\tilde\gamma}$ of restricted column-strict plane partitions
invariant under $\tilde\gamma$
and the set $\RDPP{2n-1}$
of restricted domino plane partitions with all rows of even length.
At the end of this section we state a conjecture that 
the set $\GCSPP{n,m}$ of twisted domino plane partitions
and
the set $\CDPP{n,m}$ of restricted domino plane partitions 
with all columns of even length
would have the same cardinality,
and the statistics $\overline U_{1}$ would have the same distribution.
Now we start from the definition of these classes.

\begin{definition}
Let $m$ and $n\geq1$ be nonnegative integers.
Let $\DPP{n,m}$ denote the set
of column-strict domino plane partitions $d=(d_{ij})_{1\leq i, j}$
subject to the constraints that
\begin{enumerate}
\item[(D1)]
$d$ has at most $n$ columns;
\item[(D2)]
each part in the $j$th column does not exceed $\lceil(n+m-j)/2\rceil$;
\end{enumerate}
An element of $\DPP{n,m}$ is called 
a \defterm{restricted domino plane partition} 
(abbreviated to RDPP).
If a part in the $j$th column of $c$ is equal to $\lceil(n+m-j)/2\rceil$,
we call the part \defterm{saturated}.
For $d\in\DPP{n,m}$ and a positive integer $r\geq1$,
let $\overline U_{r}(d)$ denote the number of parts equal to $r$
plus the number of saturated parts less than $r$.
Further,
if $d$ in $\DPP{n,m}$ satisfy the condition that
\begin{enumerate}
\item[(D3)]
each row (resp. column) of $d$ has even length,
\end{enumerate}
then $d$ is called a restricted column-strict domino plane partition
\defterm{with all rows (resp. columns) of even length}.
The set of all  $d\in\DPP{n,m}$
with all rows (resp. columns) of even length
is denoted by $\RDPP{n,m}$ (resp. $\CDPP{n,m}$). 
When $m=0$,
we write $\DPP{n}$ for $\DPP{n,0}$,
$\RDPP{n}$ for $\RDPP{n,0}$
and $\CDPP{n}$ for $\CDPP{n,0}$.
\end{definition}
For example,
$\RDPP{1}=\RDPP{2}=\{\emptyset\}$,
$\RDPP{3}$ is composed of the following 3 elements:
\begin{equation*}
\hbox to 20pt{\vbox{$\emptyset$,\vskip12pt}}
\qquad\qquad
\vbox{\smdomino{\h{1}\,,\\}\vskip10pt}
\qquad\qquad
\smdomino{\V{1}&\V{1}\\}\,.
\end{equation*}
$\RDPP{4}$ is composed of the following 4 elements:
\begin{equation*}
\hbox to 20pt{\vbox{$\emptyset$,\vskip12pt}}
\qquad\qquad
\vbox{\smdomino{\h{1}\,,\\}\vskip10pt}
\qquad\qquad
\smdomino{\V{1}&\V{1}\\}\,,
\qquad\qquad
\smdomino{\V{2}&\V{1}\\}\,.
\end{equation*}
$\RDPP{5}$ has 26 elements,
$\RDPP{6}$ has 50 elements,
and $\RDPP{7}$ has 646 elements.

\bigbreak

Let $n$ be a positive integer,
and
assume $c=(c_{ij})_{1\leq i,j\leq n}$ is 
in $\CSPP{2n+1}^{\tilde\gamma}$,
i.e. $c$ is invariant under each $\tilde\pi_{2r-1}$ for $r\geq1$.
For each $\tilde\pi_{2r-1}$,
we use the notation ``free'', ``paired'' and ``saturated'' 
for parts equal to $2r+1$ or $2r$ as before.
We associate to $c$
a column-strict (generalized) domino plane partition $d$ as follows.
First of all,
since $c$ is invariant under $\tilde\pi_{1}$,
which means the shape of $c_{\geq2}$ must be even.
Thus 
we remove all $1$'s from $c$
so that each row of the resulting plane partition $c_{\geq2}$ has even length.
Next, since $c_{\geq2}$ is invariant under $\tilde\pi_{3}$,
we replace paired $3$ and $2$ by a vertical domino containing $1$,
and, if row $i$ contains $k$ free $3$'s followed by $k$ free $2$'s,
then we replace them by $k$ horizontal dominoes containing $1$.
If there exists a saturated part equal to $2$,
it should be replaced by a single box containing $1$.
Next 
we replace paired $5$ and $4$ by a vertical domino containing $2$,
and, if row $i$ contains $k$ free $5$'s followed by $k$ free $4$'s,
then we replace them with $k$ horizontal dominoes containing $2$.
If there is a saturated part equal to $4$,
it is replaced by a single box containing $2$.
We repeat this process and finally obtain
a column-strict (generalized) domino plane partition $d$ with all rows of even length.
Let us denote the resulting (generalized) domino plane partition $d$ by $\DTM{2n+1}(c)$.
For example,
the left-below plane partition $c$ in $\CSPP{11}^{\tilde\gamma}$
is mapped to the right-below generalized domino plane partition $d$ in $\GCSPP{9}$
by $\DTM{11}$:
\begin{equation*}
\hbox{\vbox{\hbox to 22pt{$c=$ }\vskip65pt}}
\young{
7&7&6&6&3&2&1&1\cr
5&5&4&3&1\cr
4&3&2&2\cr
1&1\cr
}
\qquad
\hbox{\vbox{\hbox to 15pt{$d=$ }\vskip65pt}}
\hbox{
\vbox{
\hbox{
\vbox{\Horizontal{$3$}
}\kern-\futosa
\vbox{\Horizontal{$3$}
}\kern-\futosa
\vbox{\Horizontal{$1$}
}
}\kern-\futosa
\hbox{
\vbox{
\Vertical{$2$}
}\kern-\futosa
\vbox{
\Horizontal{$2$}
\kern-\futosa
\Horizontal{$1$}
}\kern-\futosa
\vbox{
\Vertical{$1$}
}}\vskip19pt
}
}
\end{equation*}
\begin{theorem}
\label{thm:bij_domino}
Let $n$ be a positive integer.
Let $c$ be a plane partition in $\CSPP{2n+1}^{\tilde\gamma}$.
Then,
 $\DTM{2n+1}(c)$
has no saturated part,
i.e. $\DTM{2n+1}(c)\in\RDPP{2n-1}$.
Thus $\DTM{2n+1}$ gives a bijection of
$\CSPP{2n+1}^{\tilde\gamma}$
onto
$\RDPP{2n-1}$.
Further we have
$\overline U_{1}(\DTM{2n+1}(c))=\overline U_{2}(c)$.
\end{theorem}
\begin{demo}{Proof}
Let $c$ be a plane partition in $\CSPP{2n+1}^{\tilde\gamma}$
and let $d=\DTM{2n+1}(c)$.
Note that a saturated part $n+1-i$ 
($i=1,\dots,n$) in $d$
 can appear in the cell $(1,2i-1)$ if there exists.
Let $\lambda$ be the shape of $c_{\geq2}$.
As a tableau can be expressed by a sequence of partitions
(see \cite[p.5]{M}),
we can write the generalized domino plane partition $d$
by a sequence
\[
\emptyset=\lambda^{(0)}\subseteq\lambda^{(1)}\subseteq\dots\subseteq\lambda^{(n)}=\lambda
\]
of partitions
in which $\lambda^{(i)}/\lambda^{(i-1)}$ consists of all the cells and dominoes 
that contain the letter $n+1-i$.
If necessary,
we add a zero at the end of $\lambda$
and we write 
$\lambda=(\lambda_{1},\lambda_{2},\dots,\lambda_{2r})$
as a sequence of even length.
Then we can write 
$
\lambda^{(i)}=(\lambda^{(i)}_{1},\dots,\lambda^{(i)}_{2r})
$
with
$
\lambda^{(i)}_{1}\geq\dots\geq\lambda^{(i)}_{2r}\geq0
$
for $i=1,\dots,n$.
For example,
the above domino tableau $d$ in $\GCSPP{9}$
is expressed by the sequence:
$\lambda^{(0)}=(0,0,0,0)$,
$\lambda^{(1)}=(4,0,0,0)$,
$\lambda^{(2)}=(4,3,1,0)$
and
$\lambda^{(3)}=(6,4,4,0)$.
If we put
$
\ell^{(i)}_{j}=\lambda^{(i)}_{j}+2r-j
$
for $i=1,\dots,n$ and $j=1,\dots,2r$,
then we have
$
\ell^{(i)}_{1}>\dots>\ell^{(i)}_{2r}\geq0
$
for $i=1,\dots,n$.
In the above example,
we have 
$\ell^{(0)}=(3,2,1,0)$,
$\ell^{(1)}=(7,2,1,0)$,
$\ell^{(2)}=(7,5,2,0)$
and
$\ell^{(3)}=(9,6,5,0)$.
When $i=0$,
we have $\ell^{(0)}_{j}=2r-j$ for $j=1,\dots,2n$,
and $\ell^{(0)}$ consist of $r$ even integers and $r$ odd integers.
For each $i=1,\dots,n$,
if there is no saturated part equal to $n+1-i$,
then $\lambda^{(i)}/\lambda^{(i-1)}$ contains only dominoes,
thus the cardinalities of odd integers and even integers in $\ell^{(i)}$
are the same as those of $\ell^{(i-1)}$.
Assume there was certain $i$ such that
$\lambda^{(i)}/\lambda^{(i-1)}$ contain a saturated part $n+1-i$ 
in $(1,2i-1)$.
Then the cardinality of odd integer would decrease by $1$,
and the cardinality of even integer would increase by $1$
if we compare $\ell^{(i)}$ with $\ell^{(i-1)}$.
Thus, if there were saturated parts,
finally $\ell^{(n)}$ would contain less odd integers than even integers.
Since $\lambda$ is even partition,
we have $\ell^{(n)}_{j}+j$ must be even for $j=1,\dots,2r$,
which implies we must have the same number of even integers and odd integers
in $\ell^{(n)}$.
This is a contradiction.
Thus we conclude that $d$ has no saturated part
and it is easy to see that $d$ is in $\RDPP{2n-1}$.
The construction of $d\in\RDPP{2n-1}$ from $c\in\CSPP{2n+1}^{\tilde\gamma}$
is clearly reversible
and give a bijection of $\CSPP{2n+1}^{\tilde\gamma}$
onto $\RDPP{2n-1}$.
Since
$3$'s and $2$'s in $c$ is replaced by dominoes containing $1$'s,
we have $\overline U_{2}(c)=\overline U_{1}(d)$.
This completes the proof.
\end{demo}

\bigbreak

Next we give examples of $\CDPP{k}$.
We have
$\CDPP{1}=\{\emptyset\}$,
$\CDPP{2}=\left\{\emptyset,\,\vbox to 18pt{\smdomino{\V{\pmb{\it1}}\\}}\,\right\}$,
and $\CDPP{3}$ has the following 3 elements:
\begin{equation*}
\hbox to 20pt{\vbox{$\emptyset$,\vskip12pt}}
\qquad\qquad
\smdomino{\V{\pmb{\it1}}\\}\,,
\qquad\qquad
\smdomino{\V{\pmb{\it1}}&\V{\pmb{\it1}}\\}\,.
\end{equation*}
$\CDPP{4}$ has the following 10 elements:
\begin{equation*}
\hbox to 20pt{\vbox{$\emptyset$,\vskip12pt}}
\qquad\qquad
\smdomino{\V{1}\,,\\}
\qquad\qquad
\smdomino{\V{\pmb{\it2}}\,,\\}
\qquad\qquad
\smdomino{\V{1}&\V{\pmb{\it1}}\,,\\}
\qquad\qquad
\smdomino{\V{\pmb{\it2}}&\V{\pmb{\it1}}\,,\\}
\end{equation*}
\begin{equation*}
\vbox{\smdomino{\V{1}&\V{\pmb{\it1}}&\V{\pmb{\it1}}\,,\\}\vskip24pt}
\qquad\qquad
\vbox{\smdomino{\V{\pmb{\it2}}&\V{\pmb{\it1}}&\V{\pmb{\it1}}\,,\\}\vskip24pt}
\qquad\qquad
\smdomino{\V{\pmb{\it2}}\\\V{1}\,,\\}
\qquad\qquad
\smdomino{\V{\pmb{\it2}}&\V{\pmb{\it1}}\\\V{1}\\},
\qquad\qquad
\smdomino{\V{\pmb{\it2}}&\V{\pmb{\it1}}&\V{\pmb{\it1}}\\\V{1}\\},
\end{equation*}
and it is not hard to see that
$\CDPP{5}$ has 25 elements,
$\CDPP{6}$ has 140 elements,
and $\CDPP{7}$ has 588 elements.
In this example the bold-faced parts are saturated.
The author computed examples for small $n$, $m$
and observe that
the cardinalities of $\GCSPP{n,m}$ and $\CDPP{n,m}$ agree 
for $n\leq6$.
\begin{conjecture}
\label{conj:bijection}
Let $m$ and $n\geq1$ be non-negative integers.
Then there would be a bijection which proves that
$\GCSPP{n,m}$
and
$\CDPP{n,m}$ has the same cardinality.
Moreover this bijection keeps $\overline U_{1}$ (i.e. the number of $1$s) invariant.
\end{conjecture}

\section{
Determinantal formulae
}\label{sec:paired}

A standard map which associate a $k$-tuple of tableaux with
a given $k$-rim hook tableaux is presented in the paper \cite[Section~6]{SW}
(see also \cite[Theorem~6.3]{CL}).
Essentially we can use this map to associate a pair of column-strict plane partitions
in $\CSPP{n,m}$ with a domino plane partition in $\DPP{n,m}$.
By this map
we can rewrite the statistics $\overline U_{k}$ on $\DPP{n,m}$ defined in the previous section
as the sum of the statistics of each column-strict plane partition.
Thus we can obtain the generating functions of $\RDPP{n,m}$ and $\CDPP{n,m}$
as an application of \cite[Lemma~7.1]{I2}
(see Theorem~\ref{thm:SW}, Corollary~\ref{cor:main_gen}).
As a corollary of these generating functions
we obtain a determinantal formula (Corollary~\ref{cor:result}(ii)) for Conjecture~\ref{conj:VS}.
Thus Conjecture~\ref{conj:VS} reduce to a determinant evaluation problem
(Conjecture~\ref{conj:result}(i)).
In the special case where $t=1$,
we prove the conjecture from Andrews-Burge determinant 
(see Lemma~\ref{lem:AB}, Theorem~\ref{thm:result}).

Let $\lambda$ be a partition.
We define a pair $(\lambda^{(0)},\lambda^{(1)})$
which is called \defterm{$2$-quotient} of $\lambda$ as follows.
If necessary,
we add a zero at the end of $\lambda$ and we regard
$\lambda=(\lambda_{1},\lambda_{2},\dots,\lambda_{2n})$
as a sequence of integers of even length.
Let 
$\ell=\lambda+\delta_{2n}
=(\lambda_{1}+2n-1,\lambda_{2}+2n-2,\dots,\lambda_{2n})$,
and put 
$\ell=\ell^{(0)}\uplus\ell^{(1)}$
where
 $\ell^{(0)}=\{x\in\ell:x\equiv0\;(\Mod2)\}$
and $\ell^{(1)}=\{x\in\ell:x\equiv1\;(\Mod2)\}$.
We can write 
$\ell^{(0)}=(2k^{(0)}_{1},\dots,2k^{(0)}_{r})$
and
$\ell^{(1)}=(2k^{(1)}_{1}+1,\dots,2k^{(1)}_{s}+1)$
where 
$k^{(0)}_{1}>\dots>k^{(0)}_{r}\geq0$
and
$k^{(1)}_{1}>\dots>k^{(1)}_{s}\geq0$.
The partition $\lambda^{(0)}$ (resp. $\lambda^{(1)}$)
is defined to be $(k^{(0)}_{1}-r+1,k^{(0)}_{2}-r+2,\dots,k^{(0)}_{r})$
(resp. $(k^{(1)}_{1}-s+1,k^{(1)}_{2}-s+2,\dots,k^{(1)}_{s})$).
For example,
if $\lambda=(5,5,3,1,1,1)$,
then we have $\lambda^{(0)}=(3,2,1)$
and $\lambda^{(1)}=(2)$.

\begin{definition}
Let $m$ and $n\geq1$ be nonnegative integers.
Let $\PCSPP{n,m}$ denote the set of all pairs $p=(c_{0},c_{1})$
of plane partitions such that
\begin{enumerate}
\item[(i)]
$c_{0}\in\CSPP{n_{0},m_{0}}$
where
$n_{0}=\left\lceil\frac{n}{2}\right\rceil$
and
$m_{0}=\left\lceil\frac{n+m+1}{2}\right\rceil-n_{0}$,
\item[(ii)]
$c_{1}\in\CSPP{n_{1},m_{1}}$
where
$n_{1}=\left\lfloor\frac{n}{2}\right\rfloor$
and
$m_{1}=\left\lfloor\frac{n+m+1}{2}\right\rfloor-n_{1}$.
\end{enumerate}
For $p\in\PCSPP{n,m}$ and a positive integer $r\geq1$,
let $\overline U_{r}(p)=\overline U_{r}(c_{0})+\overline U_{r}(c_{1})$.
Further
let $\VPCSPP{n,m}$
denote the set of all pairs $p=(c_{0},c_{1})\in\PCSPP{n,m}$
such that $\SHAPE{c_{1}}\subseteq\SHAPE{c_{0}}$ 
and $\SHAPE{c_{0}}\backslash\SHAPE{c_{1}}$
is a vertical strip.
Meanwhile,
let $\HPCSPP{n,m}$
denote the set of all pairs $p=(c_{0},c_{1})\in\PCSPP{n,m}$
such that $\SHAPE{c_{0}}\subseteq\SHAPE{c_{1}}$
and $\SHAPE{c_{1}}\backslash\SHAPE{c_{0}}$ is a horizontal strip.
We also write $\VPCSPP{n}$ for $\VPCSPP{n,0}$,
and $\HPCSPP{n}$ for $\HPCSPP{n,0}$ in short.
\end{definition}
%
%
%
For example,
if $n=4$ and $m=0$,
then
the pairs $(c_{0},c_{1})$ in $\VPCSPP{4}$
satisfy the condition that
$c_{0}\in\CSPP{2,1}$, $c_{1}\in\CSPP{2}$
and $\SHAPE{c_{0}}\backslash\SHAPE{c_{1}}$ is a vertical strip.
Thus,
$\VPCSPP{4}$ has the following 10 elements:
\begin{align*}
&\left(
\emptyset,\,
\emptyset
\right),
\qquad
\left(
\hbox{\kern2pt\vbox{\smyoung{1\\}\kern-2pt}\kern2pt},\,
\emptyset
\right),
\qquad
\left(
\hbox{\kern2pt\vbox{\smyoung{1\\}\kern-2pt}\kern2pt},\,
\hbox{\vbox{\smyoung{\pmb{\it1}\\}\kern-2pt}\kern2pt}
\right),
\qquad
\left(
\hbox{\kern2pt\vbox{\smyoung{1&\pmb{\it1}\\}\kern-2pt}\kern2pt},\,
\hbox{\vbox{\smyoung{\pmb{\it1}\\}\kern-2pt}\kern2pt}
\right),
\qquad
\left(
\hbox{\kern2pt\vbox{\smyoung{\pmb{\it2}\\}\kern-2pt}\kern2pt},\,
\emptyset
\right),
\\
&\left(
\hbox{\kern2pt\vbox{\smyoung{\pmb{\it2}\\}\kern-2pt}\kern2pt},\,
\hbox{\vbox{\smyoung{\pmb{\it1}\\}\kern-2pt}\kern2pt}
\right),
\qquad
\left(
\hbox{\kern2pt\vbox{\smyoung{\pmb{\it2}&\pmb{\it1}\\}\kern-2pt}\kern2pt},\,
\hbox{\vbox{\smyoung{\pmb{\it1}\\}\kern-2pt}\kern2pt}
\right),
\qquad
\left(
\hbox{\kern3pt\vbox{\smyoung{\pmb{\it2}\\1\\}\kern-12pt}\kern2pt},\,
\emptyset
\right),
\qquad
\left(
\hbox{\kern3pt\vbox{\smyoung{\pmb{\it2}\\1\\}\kern-13pt}\kern2pt},\,
\hbox{\vbox{\smyoung{\pmb{\it1}\\}\kern-2pt}\kern2pt}
\right),
\qquad
\left(
\hbox{\kern3pt\vbox{\smyoung{\pmb{\it2}&\pmb{\it1}\\1\\}\kern-13pt}\kern2pt},\,
\hbox{\vbox{\smyoung{\pmb{\it1}\\}\kern-2pt}\kern2pt}
\right).
\end{align*}
The italic characters stand for saturated parts.
Meanwhile,
the pairs $(c_{0},c_{1})$ in $\HPCSPP{4}$
satisfy the condition that
$c_{0}\in\CSPP{2,1}$, $c_{1}\in\CSPP{2}$
and $\SHAPE{c_{1}}\backslash\SHAPE{c_{0}}$ is a horizontal strip.
Thus,
$\HPCSPP{4}$ has the following 4 elements:
\begin{align*}
&\left(
\emptyset,\,
\emptyset
\right),
\qquad
\left(
\emptyset,\,
\hbox{\vbox{\smyoung{\pmb{\it1}\\}\kern-2pt}\kern2pt}
\right),
\qquad
\left(
\hbox{\kern2pt\vbox{\smyoung{1\\}\kern-2pt}\kern2pt},\,
\hbox{\vbox{\smyoung{\pmb{\it1}\\}\kern-2pt}\kern2pt}
\right),
\qquad
\left(
\hbox{\kern2pt\vbox{\smyoung{\pmb{\it2}\\}\kern-2pt}\kern2pt},\,
\hbox{\vbox{\smyoung{\pmb{\it1}\\}\kern-2pt}\kern2pt}
\right).
\end{align*}
%
%
%
Here we describe a bijection of column-strict domino plane partitions
onto pairs of column-strict plane partitions.
Given a domino $\alpha$,
the upper-rightmost cell in $\alpha$ is called the \newterm{head} of $\alpha$,
and
the lower-leftmost cell in $\alpha$ is called the \newterm{tail} of $\alpha$.
For any skew-shape,
we number the diagonals,
beginning with $0$ for the main diagonal,
increasing up and to the right
and decreasing down and to the left.
The \newterm{diagonal} $\diag{\alpha}$ of a domino $\alpha$ 
is the diagonal of the head of $\alpha$,
and the \newterm{color} $\Col{\alpha}$ of $\alpha$ is $\diag{\alpha}\,\Mod\,2$.
In the following pictures,
we draw only the diagonals whose numbers are even.
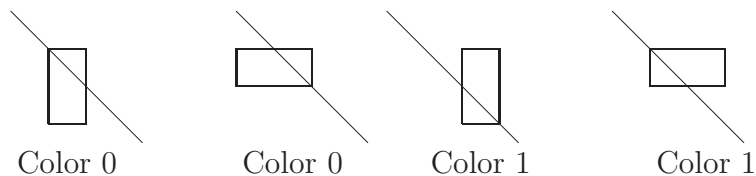
\setlength{\unitlength}{1pt} 
\setlength{\unitlength}{0.5mm}
\[
\begin{picture}(200,40)
%
%
\put( 10, 10){\line(0,1){20}}
\put( 10, 10){\line(1,0){10}}
\put( 20, 10){\line(0,1){20}}
\put( 10, 30){\line(1,0){10}}
\put(  0, 40){\line(1,-1){35}}
\put( 15,  0){\makebox(0,0){Color $0$}}
%
%
\put( 60, 20){\line(0,1){10}}
\put( 60, 20){\line(1,0){20}}
\put( 80, 20){\line(0,1){10}}
\put( 60, 30){\line(1,0){20}}
\put( 60, 40){\line(1,-1){35}}
\put( 75,  0){\makebox(0,0){Color $0$}}
%
%
\put(120, 10){\line(0,1){20}}
\put(120, 10){\line(1,0){10}}
\put(130, 10){\line(0,1){20}}
\put(120, 30){\line(1,0){10}}
\put(100, 40){\line(1,-1){35}}
\put(125,  0){\makebox(0,0){Color $1$}}
%
%
\put(170, 20){\line(0,1){10}}
\put(170, 20){\line(1,0){20}}
\put(190, 20){\line(0,1){10}}
\put(170, 30){\line(1,0){20}}
\put(160, 40){\line(1,-1){35}}
\put(185,  0){\makebox(0,0){Color $1$}}
\end{picture}
\]
%
%
%
Assume we are given a column-strict domino plane partition $d$.
For each $k=0,1$,
we associate a column-strict plane partition $c_{k}$ with $d$.
Along each diagonal $\diag{r}$,
we read only the numbers in the dominoes
that cross $\diag{r}$ and has color $k$.
We replace the dominoes by single cells 
containing the numbers and slide down along the diagonal $\diag{r}$.
In this way we obtain an (ordinary) column-strict partition $c_{k}$.
For example,
we associate the column-strict domino plane partition
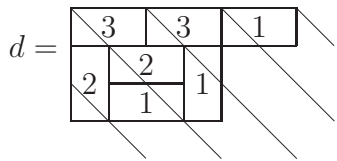
\setlength{\unitlength}{1pt} 
\setlength{\unitlength}{0.5mm}
\[
\begin{picture}(100,45)
%
%
\put( 20, 40){\line(1,0){60}}
\put( 20, 30){\line(1,0){60}}
\put( 30, 20){\line(1,0){20}}
\put( 20, 10){\line(1,0){40}}
\put( 30, 20){\line(1,0){10}}
\put( 20, 10){\line(0,1){30}}
\put( 30, 10){\line(0,1){20}}
\put( 40, 30){\line(0,1){10}}
\put( 50, 10){\line(0,1){20}}
\put( 60, 10){\line(0,1){30}}
\put( 80, 30){\line(0,1){10}}
\put( 30, 35){\makebox(0,0){$3$}}
\put( 50, 35){\makebox(0,0){$3$}}
\put( 70, 35){\makebox(0,0){$1$}}
\put( 25, 20){\makebox(0,0){$2$}}
\put( 40, 25){\makebox(0,0){$2$}}
\put( 40, 15){\makebox(0,0){$1$}}
\put( 55, 20){\makebox(0,0){$1$}}
\put( 80, 40){\line(1,-1){10}}
\put( 60, 40){\line(1,-1){30}}
\put( 40, 40){\line(1,-1){40}}
\put( 20, 40){\line(1,-1){40}}
\put( 20, 20){\line(1,-1){20}}
\put( 10, 30){\makebox(0,0){$d=$}}
\end{picture}
\]
the pair
\[
\hbox{\vbox{\hbox to 25pt{$c_{0}=$ }\vskip27pt}}
\hbox{\vbox{\young{1&1\\}\kern20pt}}
\qquad\qquad
\hbox{\vbox{\hbox to 26pt{$c_{1}=$ }\vskip27pt}}
\young{3&3&1\\2&2\\}
\]
of plane partitions.
Let $\Phi$ denote the map which associate the pair $(c_{0},c_{1})$
of column-strict plane partitions
with a column-strict domino plane partition $d$
(cf. \cite[Section~6]{SW}, \cite[Theorem~6.3]{CL}).
%
%
%
\begin{theorem}
\label{thm:SW}
The restriction of the map $\Phi:d\mapsto(c_{0},c_{1})$ to $\DPP{n,m}$
gives a bijection of
$\DPP{n,m}$ onto $\PCSPP{n,m}$.
By this bijection the $(\SHAPE{c_{0}},\SHAPE{c_{1}})$
is the $2$-quotient of $\SHAPE{d}$,
and
we have $\overline U_{r}(d)=\overline U_{r}(\Phi(d))$
for $d\in\DPP{n,m}$
and $r\geq1$.
Especially,
the restriction of the map $\Phi$ to $\RDPP{n,m}$
(resp. $\CDPP{n,m}$)
gives a bijection of
$\RDPP{n,m}$ (resp. $\CDPP{n,m}$) onto $\HPCSPP{n,m}$ (resp. $\VPCSPP{n,m}$)
which preserve the statistics $\overline U_{r}$.
\end{theorem}
\begin{demo}{Proof}
The first half of the theorem is an easy consequence of the definition.
The latter half of the theorem follows from the following proposition.
\end{demo}
\begin{prop}
\label{lem:RC}
Let $d$ be a column-strict domino plane partition,
and let $(c_{0},c_{1})=\Phi(d)$.
Then
\begin{enumerate}
\item[(i)]
All rows of $d$ have even length if,
and only if,
$\SHAPE{c_{0}}\subseteq\SHAPE{c_{1}}$
and $\SHAPE{c_{1}}\setminus\SHAPE{c_{0}}$ is a horizontal strip.
\item[(ii)]
All columns of $d$ have even length if,
and only if,
$\SHAPE{c_{1}}\subseteq\SHAPE{c_{0}}$
and $\SHAPE{c_{0}}\setminus\SHAPE{c_{1}}$ is a vertical strip.
\item[(iii)]
All rows and columns of $d$ have even length if,
and only if,
$\SHAPE{c_{0}}=\SHAPE{c_{1}}$.
\end{enumerate}
\end{prop}
\begin{demo}{Proof}
Let $\lambda=\SHAPE{d}$
and let $(\lambda^{(0)},\lambda^{(1)})$ denote the $2$-quotient
of $\lambda$.
We regard
$\lambda=(\lambda_{1},\lambda_{2},\dots,\lambda_{2n})$
as a sequence of integers of even length as before.
Put $\ell=\lambda+\delta_{2n}$,
and let $\ell^{(0)}=\{x\in\ell\,:\,x\equiv0\;(\Mod2)\}$
and
$\ell^{(1)}=\{x\in\ell\,:\,x\equiv1\;(\Mod2)\}$.
\begin{enumerate}
\item[(i)]
Note that $\lambda$ is even partition if,
and only if,
$\ell_{i}\equiv i\;(\Mod2)$.
This implies $\ell^{(0)}=(\ell_{2},\dots,\ell_{2n})=(2k^{(0)}_{1},\dots,2k^{(0)}_{n})$ and 
$\ell^{(1)}=(\ell_{1},\dots,\ell_{2n-1})=(2k^{(1)}_{1}+1,\dots,2k^{(1)}_{n}+1)$ 
with $k^{(1)}_{1}\geq k^{(0)}_{1}>k^{(1)}_{2}\geq k^{(0)}_{2}>\cdots>
k^{(1)}_{n}\geq k^{(0)}_{n}\geq0$.
Thus, if we put 
$\lambda^{(0)}=(\lambda^{(0)}_{1},\dots,\lambda^{(0)}_{n})$
and
$\lambda^{(1)}=(\lambda^{(1)}_{1},\dots,\lambda^{(1)}_{n})$,
then we have
$\lambda^{(1)}_{1}\geq\lambda^{(0)}_{1}\geq\lambda^{(1)}_{2}\geq
\lambda^{(0)}_{2}\geq\cdots\geq\lambda^{(1)}_{n}\geq\lambda^{(0)}_{n}\geq0$.
This proves that $\lambda^{(0)}\subseteq\lambda^{(1)}$
and $\lambda^{(1)}\setminus\lambda^{(0)}$ is a horizontal strip.
The reverse can be proved similarly.
\item[(ii)]
Note that $\lambda'$ is even if,
and only if,
$\ell_{2i-1}=\ell_{2i}+1$
for $i=1,\dots,n$.
This implies $\ell^{(0)}=(2k^{(0)}_{1},\dots,2k^{(0)}_{n})$ and 
$\ell^{(1)}=(2k^{(1)}_{1}+1,\dots,2k^{(1)}_{n}+1)$ 
with $k^{(0)}_{1}>\cdots>k^{(0)}_{n}\geq0$,
$k^{(1)}_{1}>\cdots>k^{(1)}_{n}\geq0$
and $k^{(1)}_{i}=k^{(0)}_{i}\text{ or }k^{(0)}_{i}-1$
for $i=1,\dots,n$.
Thus, if we put 
$\lambda^{(0)}=(\lambda^{(0)}_{1},\dots,\lambda^{(0)}_{n})$
and
$\lambda^{(1)}=(\lambda^{(1)}_{1},\dots,\lambda^{(1)}_{n})$,
we have
$\lambda^{(1)}_{i}=\lambda^{(0)}_{i}\text{ or }\lambda^{(0)}_{i}-1$
which is equivalent that $\lambda^{(1)}\subseteq\lambda^{(0)}$
and $\lambda^{(0)}\setminus\lambda^{(1)}$ is a vertical strip.
\item[(iii)]
This immediately follows from (i) and (ii).
\end{enumerate}
This completes the proof.
\end{demo}
\begin{lemma}
\label{lem:matrix}
Let $n$ and $N$ be positive integers such that
$n\leq N$.
Let $A=(a_{ij})$ and $B=(b_{ij})$
be $n\times N$ matrices.
\begin{enumerate}
\item[(i)]\label{roster:lem1}
Let $C=(c_{ij})$ denote the $n\times N$ matrix
whose $(i,j)$th entry is
$
c_{ij}=\sum_{k=1}^{j}b_{ik}.
$
Then we have
\begin{equation}
\label{eq:cauchy1}
\sum
\det\begin{pmatrix}
a_{1,j_{1}}&\hdots&a_{1,j_{n}}\\
   \vdots  &\ddots&   \vdots  \\
a_{n,j_{1}}&\hdots&a_{n,j_{n}}\\
\end{pmatrix}
\det\begin{pmatrix}
b_{1,k_{1}}&\hdots&b_{1,k_{n}}\\
   \vdots  &\ddots&   \vdots  \\
b_{n,k_{1}}&\hdots&b_{n,k_{n}}\\
\end{pmatrix}
=\det A\,{}^t\!C,
\end{equation}
where the sum runs over all $(j_{1},\dots,j_{n})$
and $(k_{1},\dots,k_{n})$ such that
\begin{equation*}
1\leq k_{1}\leq j_{1}<k_{2}\leq j_{2}<\cdots<k_{n}\leq j_{n}\leq N.
\end{equation*}
\item[(ii)]\label{roster:lem2}
Let $C'=(c_{ij}')$ denote the $n\times N$ matrix
whose $(i,j)$th entry is $c_{ij}'=\sum_{k=\max(1,j-1)}^{j}b_{ik}$.
Then we have
\begin{equation}
\label{eq:cauchy2}
\sum
\det\begin{pmatrix}
a_{1,j_{1}}&\hdots&a_{1,j_{n}}\\
   \vdots  &\ddots&   \vdots  \\
a_{n,j_{1}}&\hdots&a_{n,j_{n}}\\
\end{pmatrix}
\det\begin{pmatrix}
b_{1,k_{1}}&\hdots&b_{1,k_{n}}\\
   \vdots  &\ddots&   \vdots  \\
b_{n,k_{1}}&\hdots&b_{n,k_{n}}\\
\end{pmatrix}
=\det A\,{}^t\!C',
\end{equation}
where the sum runs over all $(j_{1},\dots,j_{n})$
and $(k_{1},\dots,k_{n})$ such that
\begin{equation*}
1\leq j_{1}<j_{2}<\cdots<j_{n}\leq N,
\end{equation*}
and $\max(j_{\nu-1}+1,j_{\nu}-1)\leq k_{\nu}\leq j_{\nu}$
for $\mu=1,2,\dots,n$.
Here we use the convention that $j_{0}=1$.
\end{enumerate}
\end{lemma}
\begin{demo}{Proof}
The both of the identities \thetag{\ref{eq:cauchy1}}, \thetag{\ref{eq:cauchy2}}
reduce to the following Cauchy-Binet identity:
Let $X=(x_{ij})$ and $Y=(y_{ij})$
be $n\times N$ matrices.
Then
\begin{equation}
\label{eq:cauchy}
\sum
\det\begin{pmatrix}
x_{1,j_{1}}&\hdots&x_{1,j_{n}}\\
   \vdots  &\ddots&   \vdots  \\
x_{n,j_{1}}&\hdots&x_{n,j_{n}}\\
\end{pmatrix}
\det\begin{pmatrix}
y_{1,j_{1}}&\hdots&y_{1,j_{n}}\\
   \vdots  &\ddots&   \vdots  \\
y_{n,j_{1}}&\hdots&y_{n,j_{n}}\\
\end{pmatrix}
=\det X\,{}^t\!Y,
\end{equation}
where the sum runs over all 
$
1\leq j_{1}<j_{2}<\cdots<j_{n}\leq N.
$
For example,
to prove \thetag{\ref{eq:cauchy1}} using \thetag{\ref{eq:cauchy}},
it suffices to show that
\[
\det C_{j_{1},j_{2},\dots,j_{n}}
=\sum_{k_{1}=1}^{j_{1}}\sum_{k_{2}=j_{1}+1}^{j_{2}}\cdots\sum_{k_{n}=j_{n-1}+1}^{j_{n}}
B_{k_{1},k_{2},\dots,k_{n}}.
\]
This can be easily seen
if one writes
\[
c_{i,j_{\nu}}
=\sum_{k=1}^{j_{\nu}}b_{i,k}
=\sum_{k=1}^{j_{1}}b_{i,k}+\sum_{k=j_{1}+1}^{j_{2}}b_{i,k}
+\cdots+\sum_{k=j_{\nu-1}+1}^{j_{\nu}}b_{i,k}
\]
for $1\leq i,\nu\leq n$.
The other identity \thetag{\ref{eq:cauchy2}} can be derived similarly.
\end{demo}
In \cite{I2} we obtained the generating function of $\CSPP{n,m}$ weighted by 
the statistics $\overline U_{k}$ and parts.
Here we cite the theorem as the following lemma without proof.
\begin{lemma}
\label{lem:lp}
(\cite[Lemma~7.1]{I2})
Let $m$ and $n\geq1$ be non-negative integers,
and fix a positive integer $N\geq n+m$.
Let $\lambda$ be a partition with $\ell(\lambda)\leq n$.
For $c\in\CSPP{n,m}$,
let $\boldsymbol{t}^{\overline U(c)}\boldsymbol{x}^c$ denote 
$
\prod_{k=1}^{N}t_{k}^{\overline U_{k}(c)}
\prod_{i\geq1}x_{i}^{m_{i}}
$,
where $m_{i}$ denote the number of $i$'s in $c$.
We put $y_{i}=t_{i}x_{i}$ and $Y_{i}=\prod_{k=i}^{N}t_{k}\,x_{i}$,
and use the vector notation $\boldsymbol{y}^{(r)}$
for the $r$-tuples $(y_{1},\dots,y_{r})$ of the variables.
Then
the generating function of all plane partitions $c\in\CSPP{n,m}$
of shape $\lambda'$ with the weight $\boldsymbol{t}^{\overline U(c)}\boldsymbol{x}^c$ is given by
\begin{equation}
\sum_{{c\in\CSPP{n,m}}\atop{\SHAPE{c}=\lambda'}}\boldsymbol{t}^{\overline U(c)}\boldsymbol{x}^c
=\det\left(e^{(n+m-i)}_{\lambda_{j}-j+i}(\boldsymbol{y}^{(n+m-i-1)},Y_{n+m-i})\right)_{1\leq i,j\leq n}.
\label{eq:GenFunc}
\end{equation}
\end{lemma}
As an application of Lemma~\ref{lem:matrix} and Lemma~\ref{lem:lp},
we obtain the following theorem,
which give us the generating function of $\RDPP{n,m}$
and $\CDPP{n,m}$ with the same weights.
%
%
\begin{theorem}
\label{thm:main_gen}
Let $m$ and $n\geq1$ be nonnegative integers,
and fix a positive integer $N\geq n+m$.
Put
$n_{0}=\left\lceil\frac{n}{2}\right\rceil$,
$m_{0}=\left\lceil\frac{n+m+1}{2}\right\rceil-n_{0}$,
$n_{1}=\left\lfloor\frac{n}{2}\right\rfloor$
and
$m_{1}=\left\lfloor\frac{n+m+1}{2}\right\rfloor-n_{1}$.
For a domino plane partition
$d\in\DPP{n,m}$,
let $\boldsymbol{t}^{\overline U(d)}\boldsymbol{x}^d$ denote 
$
\prod_{k=1}^{N}t_{k}^{\overline U_{k}(d)}
\prod_{i\geq1}x_{i}^{m_{i}}
$
where $m_{i}$ denote the number of $i$'s in $d$.
We put $y_{i}=t_{i}x_{i}$ and $Y_{i}=\prod_{k=i}^{N}t_{k}\,x_{i}$,
and write $\boldsymbol{y}^{(r)}$
for the $r$-tuples $(y_{1},\dots,y_{r})$,
where we use the convention that $(\boldsymbol{y}^{(r-1)},Y_{r})$ is empty if $r\leq0$.
\begin{enumerate}
\item[(i)]
%
%
If $n$ is even (i.e. $n_{0}=n_{1}$),
then we obtain 
$
\sum_{d\in\RDPP{n,m}}\boldsymbol{t}^{\overline U(d)}\boldsymbol{x}^{d}=\det R,
$
where 
$R=(R_{ij})_{0\leq i\leq n_{0}-1,\,0\leq j\leq n_{1}-1}$ 
is the $n_{0}\times n_{1}$ matrix
whose $(i,j)$th entry is
\begin{equation}
R_{ij}=
\sum_{k\geq0}
e^{(m_{0}+i+1)}_{k-i}(\boldsymbol{y}^{(m_{0}+i-1)},Y_{m_{0}+i},1)\,
e^{(m_{1}+j)}_{k-j}(\boldsymbol{y}^{(m_{1}+j-1)},Y_{m_{1}+j}).
\label{eq:R_even}
\end{equation}
%
%
\item[(ii)]
If $n$ is odd (i.e. $n_{0}=n_{1}+1$),
then we obtain 
$
\sum_{d\in\RDPP{n,m}}\boldsymbol{t}^{\overline U(d)}\boldsymbol{x}^{d}
=\det\left(\begin{array}{c|c}\vec r&R\end{array}\right),
$
 where 
\begin{equation}
R_{ij}=
\sum_{k\geq0}
e^{(m_{0}+i+1)}_{k-i}(\boldsymbol{y}^{(m_{0}+i-1)},Y_{m_{0}+i},1)\,
e^{(m_{1}+j)}_{k-j-1}(\boldsymbol{y}^{(m_{1}+j-1)},Y_{m_{1}+j}),
\label{eq:R_odd}
\end{equation}
and $\vec r=(r_{i})_{0\leq i\leq n_{0}-1}$ is the column vector
whose $i$th entry is $r_{i}=\delta_{i,0}$.
%
%
\item[(iii)]
If $n$ is even (i.e. $n_{0}=n_{1}$),
then we obtain 
$
\sum_{d\in\CDPP{n,m}}\boldsymbol{t}^{\overline U(d)}\boldsymbol{x}^{d}=\det C,
$
where 
$C=(C_{ij})_{0\leq i\leq n_{0}-1,\,0\leq j\leq n_{1}-1}$ is the $n_{0}\times n_{1}$ matrix
whose $(i,j)$th entry is
\begin{equation}
C_{ij}=
\sum_{k\geq0}\sum_{\nu=0}^{k}
e^{(m_{0}+i)}_{k-i}(\boldsymbol{y}^{(m_{0}+i-1)},Y_{m_{0}+i})\,
e^{(m_{1}+j)}_{\nu-j}(\boldsymbol{y}^{(m_{1}+j-1)},Y_{m_{1}+j}).
\label{eq:C_even}
\end{equation}
%
%
\item[(iv)]
If $n$ is odd (i.e. $n_{0}=n_{1}+1$),
then we obtain 
$
\sum_{d\in\CDPP{n,m}}\boldsymbol{t}^{\overline U(d)}\boldsymbol{x}^{d}
=\det\left(\begin{array}{c|c}\vec c&C\end{array}\right),
$
where 
\begin{equation}
C_{ij}=
\sum_{k\geq0}\sum_{\nu=0}^{k-1}
e^{(m_{0}+i)}_{k-i}(\boldsymbol{y}^{(m_{0}+i-1)},Y_{m_{0}+i})\,
e^{(m_{1}+j)}_{\nu-j}(\boldsymbol{y}^{(m_{1}+j-1)},Y_{m_{1}+j}),
\label{eq:C_odd}
\end{equation}
and
$\vec c=(c_{i})_{0\leq i\leq n_{0}-1}$ is the column vector
whose $i$th entry is
\[
c_{i}=
\sum_{k\geq0}e^{(m_{0}+i)}_{k-i}(\boldsymbol{y}^{(m_{0}+i-1)},Y_{m_{0}+i}).
\]
\end{enumerate}
\end{theorem}
%
%
%
\begin{demo}{Proof of Theorem~\ref{thm:main_gen}}
Let $d\in\RDPP{n,m}$.
By the bijection $\Phi:d\mapsto(c_{0},c_{1})$ in Theorem~\ref{thm:SW}
and the fact that
$\boldsymbol{t}^{\overline U(d)}\boldsymbol{x}^{d}
=\boldsymbol{t}^{\overline U(c_{0})}\boldsymbol{x}^{c_{0}}\cdot
\boldsymbol{t}^{\overline U(c_{1})}\boldsymbol{x}^{c_{1}},
$
we obtain that
$
\sum_{d\in\RDPP{n,m}}\boldsymbol{t}^{\overline U(d)}\boldsymbol{x}^{d}
=\sum_{(c_{0},c_{1})\in\HPCSPP{n,m}}
\boldsymbol{t}^{\overline U(c_{0})}\boldsymbol{x}^{c_{0}}\cdot\boldsymbol{t}^{\overline U(c_{1})}\boldsymbol{x}^{c_{1}}
$.
Write $N_{0}=n_{0}+m_{0}$ and $N_{1}=n_{1}+m_{1}$ for short.
Let $\lambda^{(0)}=\SHAPE{c_{0}}'$ and $\lambda^{(1)}=\SHAPE{c_{1}}'$.
Note that $\SHAPE{c_{1}}\setminus\SHAPE{c_{0}}$ is a horizontal strip 
if, and only if, $\lambda^{(1)}\setminus\lambda^{(0)}$ is a vertical strip.
From Lemma~\ref{lem:lp} above,
we conclude that
$\sum_{d\in\RDPP{n,m}}\boldsymbol{t}^{\overline U(d)}\boldsymbol{x}^{d}$
is equal to
\begin{align*}
\sum
\det_{1\leq i,j\leq n_{0}}\left(e^{(N_{0}-i)}_{\lambda^{(0)}_{j}-j+i}(\boldsymbol{y}^{(N_{0}-i-1)},Y_{N_{0}-i})\right)
\det_{1\leq i,j\leq n_{1}}\left(e^{(N_{1}-i)}_{\lambda^{(1)}_{j}-j+i}(\boldsymbol{y}^{(N_{1}-i-1)},Y_{N_{1}-i})\right),
\end{align*}
where the sum runs over all pair $(\lambda^{(0)},\lambda^{(1)})$
of partitions such that 
$\lambda^{(0)}\subseteq\lambda^{(1)}$ and
$\lambda^{(1)}\setminus\lambda^{(0)}$ is a vertical strip.
Rearrange the row and column indices,
then this equals
\begin{align*}
\sum
\det_{0\leq i,j\leq n_{0}-1}\left(e^{(m_{0}+i)}_{\lambda^{(0)}_{n_{0}-j}+j-i}(\boldsymbol{y}^{(m_{0}+i-1)},Y_{m_{0}+i})\right)
\det_{0\leq i,j\leq n_{1}-1}\left(e^{(m_{1}+i)}_{\lambda^{(1)}_{n_{1}-j}+j-i}(\boldsymbol{y}^{(m_{1}+i-1)},Y_{m_{1}+i})\right).
\end{align*}
If $n_{0}=n_{1}$,
put $l_{j+1}^{(0)}=\lambda^{(0)}_{n_{0}-j}+j$ and $l_{j+1}^{(1)}=\lambda^{(1)}_{n_{1}-j}+j$
for $0\leq j\leq n_{0}-1$.
Note that $0\leq l^{(k)}_{1}<\cdots<l^{(k)}_{n_{0}}$ for $k=0,1$.
It is easy to see that $\lambda^{(1)}\setminus\lambda^{(0)}$ is a vertical strip
if, and only if, $l_{j}^{(1)}-1\leq l_{j}^{(0)}\leq l_{j}^{(1)}$.
Thus,
if we put
\[
P_{ij}=
e^{(m_{1}+i)}_{j-i}(\boldsymbol{y}^{(m_{1}+i-1)},Y_{m_{1}+i})
\]
and
\begin{align*}
Q_{ij}&=
e^{(m_{0}+i)}_{j-i}(\boldsymbol{y}^{(m_{0}+i-1)},Y_{m_{0}+i})
+e^{(m_{0}+i)}_{j-i-1}(\boldsymbol{y}^{(m_{0}+i-1)},Y_{m_{0}+i})
\\
&
=e^{(m_{0}+i+1)}_{j-i}(\boldsymbol{y}^{(m_{0}+i-1)},Y_{m_{0}+i},1)
\end{align*}
then,
by Lemma~\ref{lem:matrix}(ii),
we obtain 
$\sum_{d\in\RDPP{n,m}}\boldsymbol{t}^{\overline U(d)}\boldsymbol{x}^{d}
=\det(Q_{ij}){}^{t}\!(P_{ij})$,
which leads to the desired identity \thetag{\ref{eq:C_even}}.
If $n_{0}=n_{1}+1$,
then use
\[
\det\left(A\right)
=\det\left(\begin{array}{c|c}
1&O\\\hline
O&A
\end{array}\right)
\]
and repeat the same arguments to obtain \thetag{\ref{eq:C_odd}}.
The other identities \thetag{\ref{eq:R_even}} \thetag{\ref{eq:R_odd}}
can be derived similarly using Lemma~\ref{lem:matrix}(i).
\end{demo}
If we specialize the variables in Theorem~\ref{thm:main_gen},
then we obtain the following corollary:
%
%
\begin{corollary}
\label{cor:main_gen}
Let $m$ and $n\geq1$ be nonnegative integers,
and fix a positive integer $k\geq1$.
Put
$n_{0}=\left\lceil\frac{n}{2}\right\rceil$,
$m_{0}=\left\lceil\frac{n+m+1}{2}\right\rceil-n_{0}$,
$n_{1}=\left\lfloor\frac{n}{2}\right\rfloor$
and
$m_{1}=\left\lfloor\frac{n+m+1}{2}\right\rfloor-n_{1}$.
\begin{enumerate}
%
%
\item[(i)]
If $n$ is even (i.e. $n_{0}=n_{1}$),
then we obtain 
$
\sum_{d\in\RDPP{n,m}}t^{\overline U_{k}(d)}=\det R'(t),
$
where 
$R'(t)=(R_{ij}')_{0\leq i\leq n_{0}-1,\,0\leq j\leq n_{1}-1}$ is the $n_{0}\times n_{1}$ matrix
whose $(i,j)$th entry is
\begin{align}
R_{ij}'
&=
\binom{m_{0}+m_{1}+i+j-1}{m_{0}+2i-j} (1+t^2)
\nonumber\\
&\qquad
+\left\{\binom{m_{0}+m_{1}+i+j-1}{m_{0}+2i-j-1}+\binom{m_{0}+m_{1}+i+j-1}{m_{0}+2i-j+1}\right\}t
\label{eq:binomt_even}
\end{align}
if $m_{0}+i>0$ and $m_{1}+j>0$,
$
R_{0,j}'=\binom{m_{1}+j}{-j+1}+\binom{m_{1}+j}{-j}t
$
if $m_{0}=0$ and $m_{1}+j>0$,
$
R_{i,0}'=\delta_{0,i}
$
if  and $m_{0}+i\geq0$ and $m_{1}=0$.
%
%
\item[(ii)]
If $n$ is odd (i.e. $n_{0}=n_{1}+1$),
then we obtain 
$
\sum_{d\in\RDPP{n,m}}t^{\overline U_{k}(d)}
=\det\left(\begin{array}{c|c}\vec r&R'(t)\end{array}\right),
$
where 
$R'(t)=(R_{ij}')_{0\leq i\leq n_{0}-1,\,0\leq j\leq n_{1}-1}$ is the $n_{0}\times n_{1}$ matrix
whose $(i,j)$th entry is
\begin{align}
R_{ij}'
&=
\binom{m_{0}+m_{1}+i+j-1}{m_{0}+2i-j-1} (1+t^2)
\nonumber\\
&\qquad
+\left\{\binom{m_{0}+m_{1}+i+j-1}{m_{0}+2i-j-2}+\binom{m_{0}+m_{1}+i+j-1}{m_{0}+2i-j}\right\}t
\label{eq:binomt_odd}
\end{align}
if $m_{0}+i>0$ and $m_{1}+j>0$,
$
R_{i,0}'=\binom{m_{0}+i}{-i+1}+\binom{m_{0}+i}{-i}t
$
if $m_{0}+i>0$ and $m_{1}=0$,
and
$
R_{0,j}'=\delta_{0,j}
$
if $m_{0}=0$ and $m_{1}+j\geq0$.
Here $\vec r=(\delta_{i,0})_{0\leq i\leq n_{0}-1}$ is as in Theorem~\ref{thm:main_gen}.
%
%
\item[(iii)]
If $n$ is even (i.e. $n_{0}=n_{1}$),
then we obtain 
$
\sum_{d\in\CDPP{n,m}}t^{\overline U_{k}(d)}=\det C'(t),
$
where 
$C'(t)=(C_{ij}')_{0\leq i\leq n_{0}-1,\,0\leq j\leq n_{1}-1}$ is the $n_{0}\times n_{1}$ matrix
whose $(i,j)$th entry is
\begin{align}
C_{ij}'&=
\sum_{k\leq m_{0}+2i-j-1}\left[
\binom{m_{0}+m_{1}+i+j-2}{k}(1+t^2)
\right.\nonumber\\
&\left.
+\left\{\binom{m_{0}+m_{1}+i+j-2}{k-1}+\binom{m_{0}+m_{1}+i+j-2}{k+1}\right\}t
\right]
\label{eq:c_{ij}even}
\end{align}
if $m_{0}+i>0$ and $m_{1}+j>0$,
$C_{i,0}'=2^{m_{0}+i-1}(1+t)$ if $m_{0}+i>0$ and $m_{1}=0$
and
$C_{0,j}'=\delta_{0,j}$ if $m_{0}=0$ and $m_{1}+j\geq0$.
%
%
\item[(iv)]
If $n$ is odd (i.e. $n_{0}=n_{1}+1$),
then we obtain 
$
\sum_{d\in\CDPP{n,m}}t^{\overline U_{k}(d)}
=\det\left(\begin{array}{c|c}\vec c'(t)&C'(t)\end{array}\right),
$
where 
$C'(t)=(C_{ij}')_{0\leq i\leq n_{0}-1,\,0\leq j\leq n_{1}-1}$ is the $n_{0}\times n_{1}$ matrix
whose $(i,j)$th entry is
\begin{align}
C_{ij}'&=
\sum_{k<m_{0}+2i-j-1}\left[
\binom{m_{0}+m_{1}+i+j-2}{k}(1+t^2)
\right.\nonumber\\
&\left.
+\left\{\binom{m_{0}+m_{1}+i+j-2}{k-1}+\binom{m_{0}+m_{1}+i+j-2}{k+1}\right\}t
\right]
\label{eq:c_{ij}odd}
\end{align}
if $m_{0}+i>0$ and $m_{1}+j>0$,
$2^{m_{0}+i-1}(1+t)-\delta_{i,0}$ if $m_{0}+i>0$ and $m_{1}+j=0$,
and $0$ if $m_{0}+i=0$ and $m_{1}+j\geq0$.
Here 
$\vec c'(t)=(c_{i}')_{0\leq i\leq n_{0}-1}$ is the column vector
whose $i$th entry is
$
c_{i}'=\begin{cases}
1,
&\text{ if $m_{0}=i=0$,}\\
2^{m_{0}+i-1}\left(1+t\right),
&\text{ otherwise.}
\end{cases}
$
\end{enumerate}
\end{corollary}
For example,
if $n=5$ and $m=0$, then
$n_{0}=3$, $m_{0}=0$, $n_{1}=2$, $m_{1}=1$,
and $\sum_{d\in\RDPP{5}}t^{\overline U_{k}(d)}$, $k\geq1$, is equal to
\[
\det\left(
\begin {array}{ccc} 
1&1&0\\
0&1+t+t^{2}&1+2t+t^{2}\\
0&t&3+4t+3t^{2}
\end {array}
\right)
=3+6t+8t^{2}+6t^{3}+3t^{4},
\]
whereas
$\sum_{d\in\CDPP{5}}t^{\overline U_{k}(d)}$, $k\geq1$, is equal to
\[
\det\left(
\begin {array}{ccc} 
1&0&0\\
1+t&1+t+t^{2}&t\\
2+2t&2+4t+2t^{2}&3+5t+3t^{2}
\end {array}\right)
=3+6t+7t^{2}+6t^{3}+3t^{4}.
\]
%
%
%
\begin{demo}{Proof of Corollary~\ref{cor:main_gen}}
Substitute $x_{i}=1$ for $i\geq1$, $t_{k}=t$ and $t_{i}=1$ for $i\neq k$ into Theorem~\ref{thm:main_gen}.
First assume $n_{0}=n_{1}$.
Using 
$
e^{n}_{r}(t,1,\dots,1)
=\begin{cases}
\delta_{0,r},&\text{ if $n=0$,}\\
\binom{n-1}{r}+\binom{n-1}{r-1}t,
&\text{ if $n>0$,}
\end{cases}
$
we obtain
\[
R_{i,j}'=\sum_{k\geq0}
\left\{
\binom{m_{0}+i}{k-i}+\binom{m_{0}+i}{k-i-1}t
\right\}
\left\{
\binom{m_{1}+j-1}{k-j}+\binom{m_{1}+j-1}{k-i-1}t
\right\}
\]
if $m_{1}+j>0$.
Use the binomial coefficient identity 
$
\sum_{k\geq0}\binom{\alpha}{k}\binom{\beta}{\gamma-k}=\binom{\alpha+\beta}{\gamma}
$
to obtain \thetag{\ref{eq:binomt_even}}.
The case of $m_{1}+j=0$ is treated similarly. 
One can prove \thetag{\ref{eq:binomt_odd}} similarly in the case of $n_{0}=n_{1}+1$.
Next,
to prove (ii)(a),
substitute $x_{i}=1$ for $i\geq1$, $t_{k}=t$ and $t_{i}=1$ for $i\neq k$ 
into \thetag{\ref{eq:C_even}}.
If $m_{0}+i>0$ and $m_{1}+j>0$, then
we obtain
\begin{align*}
C_{i,j}'=\sum_{k\geq0}\sum_{\nu=0}^{k}
\left\{
\binom{m_{0}+i-1}{k-i}+\binom{m_{0}+i-1}{k-i-1}t
\right\}
\left\{
\binom{m_{1}+j-1}{\nu-j}+\binom{m_{1}+j-1}{\nu-i-1}t.
\right\}
\end{align*}
Thus,
the proof of the desired identity \thetag{\ref{eq:c_{ij}even}}
reduce to the following identity
\[
\sum_{k\geq0}\sum_{\nu=0}^{k}
\binom{m_{0}+i-1}{k-i-\alpha}\binom{m_{1}+j-1}{\nu-j-\beta}
=\sum_{k< m_{0}+2i-j+\alpha-\beta}\binom{m_{0}+m_{1}+i+j-2}{k}.
\]
The other cases can be treated similarly.
\end{demo}
If we put $m=0$ in Corollary~\ref{cor:main_gen},
then we obtain the following corollary:
\begin{corollary}
\label{cor:result}
Let $n$ be a positive integer,
and fix a positive integer $k\geq1$.
Put
$n_{0}=\left\lceil\frac{n}{2}\right\rceil$.
\begin{enumerate}
%
%
\item[(i)]
If $n$ is even,
then we obtain 
$
\sum_{d\in\RDPP{n}}t^{\overline U_{k}(d)}=\det R_{n_{0}}^\text{e}(t),
$
where 
$R_{n_{0}}^\text{e}(t)=(R^\text{e}_{ij})_{1\leq i,j\leq n_{0}-1}$ is 
the $n_{0}\times n_{0}$ matrix
whose $(i,j)$th entry is
\begin{equation}
R^\text{e}_{ij}=\binom{i+j}{2i-j+1}(1+t^2)
+\left\{
\binom{i+j}{2i-j}+\binom{i+j}{2i-j+2}
\right\}t
\end{equation}
if $(i,j)\neq(0,0)$ and $R^\text{e}_{0,0}=1$.
%
%
\item[(ii)]
If $n$ is odd,
then we obtain 
$
\sum_{d\in\RDPP{n}}t^{\overline U_{k}(d)}=\det R_{n_{0}}^\text{o}(t),
$
where 
$R_{n_{0}}^\text{o}(t)=(R^\text{o}_{ij})_{0\leq i,j\leq n_{0}-1}$ is 
the $n_{0}\times n_{0}$ matrix
whose $(i,j)$th entry is
\begin{equation}
R^\text{o}_{ij}=\binom{i+j-1}{2i-j}(1+t^2)
+\left\{
\binom{i+j-1}{2i-j-1}+\binom{i+j-1}{2i-j+1}
\right\}t
\end{equation}
if $(i,j)\neq(0,0),(0,1)$ and 
$R^\text{o}_{0,0}=R^\text{o}_{0,1}=1$.
%
%
\item[(iii)]
If $n$ is even,
then we obtain 
$
\sum_{d\in\CDPP{n}}t^{\overline U_{k}(d)}=\det C_{n_{0}}^\text{e}(t),
$
where 
$C_{n_{0}}^\text{e}(t)=(C_{ij}^\text{e})_{0\leq i,j\leq n_{0}-1}$ is 
the $n_{0}\times n_{0}$ matrix
whose $(i,j)$th entry is
\begin{align}
C_{ij}^\text{e}&=\left\{2\binom{i+j-2}{2i-j-1}+\binom{i+j-2}{2i-j}\right\}(1+t^2)
\\
&+\left\{
2\binom{i+j-2}{2i-j-2}+\binom{i+j-2}{2i-j-1}+2\binom{i+j-2}{2i-j}+\binom{i+j-2}{2i-j+1}
\right\}t
\nonumber
\end{align}
if $i+j\geq2$, $C^\text{e}_{0,0}=1+t$, $C^\text{e}_{0,1}=t$ and $C^\text{e}_{1,0}=0$.
\item[(iv)]
If $n$ is odd,
then we obtain 
$
\sum_{d\in\CDPP{n}}t^{\overline U_{k}(d)}=\det C_{n_{0}}^\text{o}(t),
$
where 
$C_{n_{0}}^\text{o}(t)=(C^\text{o}_{ij})_{0\leq i,j\leq n_{0}-1}$ is 
the $n_{0}\times n_{0}$ matrix
whose $(i,j)$th entry is
\begin{align}
C^\text{o}_{ij}&=\left\{2\binom{i+j-3}{2i-j-2}+\binom{i+j-3}{2i-j-1}\right\}(1+t^2)
\\
&+\left\{
2\binom{i+j-3}{2i-j-3}+\binom{i+j-3}{2i-j-2}+2\binom{i+j-3}{2i-j-1}+\binom{i+j-3}{2i-j}
\right\}t
\nonumber\end{align}
\end{enumerate}
if $i+j\geq3$, $C^\text{o}_{0,0}=1$, $C^\text{o}_{0,1}=C^\text{o}_{0,2}=C^\text{o}_{2,0}=0$, $C^\text{o}_{1,0}=1+t$ and $C^\text{o}_{1,1}=1+t+t^2$.
\end{corollary}
\begin{demo}{Proof}
Substitute $m=0$ in Corollary~\ref{cor:main_gen}.
Then we have $m_{0}=1$ and $m_{1}=0$ if $n$ is even,
and we have $m_{0}=0$ and $m_{1}=1$ if $n$ is odd.
Using this,
we can directly compute each entry of $R_{n_{0}}^\text{e}(t)=R'(t)$ 
(resp. $R_{n_{0}}^\text{o}(t)=\left(\begin{array}{c|c}\vec r&R'\end{array}(t)\right)$)
in (i)(ii) of Corollary~\ref{cor:main_gen}
if $n$ is even (resp. odd).
We define $n\times n$ matrices $U^\text{e}_{n}$ and $U^\text{o}_{n}$ by
$
U^\text{e}_{n}=\left(\delta_{i,j}-2\delta_{i+1,j}\right)_{0\leq i,j\leq n-1}
$
and
$
U^\text{o}_{n}=U^\text{e}_{n}+\left(2\delta_{i^2+(j-1)^2,0}\right)_{0\leq i,j\leq n-1}.
$
As before we substitute $m=0$ in (iii)(iv) of Corollary~\ref{cor:main_gen}.
Then put $C_{n_{0}}^\text{e}(t)=U^\text{e}_{n_{0}}C'(t)$ 
(resp. $C_{n_{0}}^\text{o}(t)=U^\text{o}_{n_{0}}\left(\begin{array}{c|c}\vec r&R'(t)\end{array}\right)$)
if $n$ is even (resp. odd).
We can compute each entry of $C_{n_{0}}^\text{e}(t)$ and $C_{n_{0}}^\text{o}(t)$ directly.
\end{demo}
For instance, 
in the case where $n=5$ and $n=6$,
$R_{3}^\text{o}$ and $R_{3}^\text{e}$ are 
\[
\left(
\begin {array}{ccc}
1&1&0\\
0&1+t+{t}^{2}&1+2t+{t}^{2}\\
0&t&3+4t+3{t}^{2}
\end {array}
\right),
\quad
\left( 
\begin {array}{ccc} 
1&1+t+{t}^{2}&t\\
0&1+2t+{t}^{2}&3+4t+3{t}^{2}\\
0&t&4+7t+4{t}^{2}
\end {array} \right),
\]
and $C^\text{o}$ and $C^\text{e}$ are as follows:
\[
\left(
\begin {array}{ccc}
1&0&0\\
1+t&1+t+{t}^{2}&t\\
0&2t&3+3t+3{t}^{2}
\end {array}\right),
\quad
\left( \begin {array}{ccc}
1+t&t&0\\
0&2+t+2{t}^{2}&1+3t+{t}^{2}\\
0&2t&5+6t+5{t}^{2}
\end {array} \right).
\]
Thus we obtain certain determinantal formulae
for the generating functions of $\RDPP{n}$ and $\CDPP{n}$
weighted by $\overline U_{k}$.
Thus the above (ii) with Theorem~\ref{cor:bijection_TSPPtoCSPP},
Theorem~\ref{thm:bij_domino} and Theorem~\ref{thm:SW}
proves Theorem~\ref{thm:VS}.
Here the problem of evaluations still remains:
\begin{conjecture}
\label{conj:result}
Let $r$ be a positive integer.
Let $R_{r}^\text{o}$,
$C_{r}^\text{e}$ and $C_{r}^\text{o}$ as in Corollary~\ref{cor:result}
Then the following identities would hold.
\begin{enumerate}
%
%
\item[(i)]
For $r\geq1$,
\begin{equation}
\det R_{r}^\text{o}(t)=A^\text{VS}_{2r+1}(t).
\label{eq:R_o}
\end{equation}
\item[(ii)]
For $r\geq1$,
\begin{equation}
\det C_{r}^\text{e}(t)=A^\text{HTS}_{2r}(t).
\label{eq:C_e}
\end{equation}
\item[(iii)]
For $r\geq1$,
\begin{equation}
\det C_{r}^\text{o}(t)=A^\text{HTS}_{2r-1}(t).
\label{eq:C_o}
\end{equation}
\end{enumerate}
\end{conjecture}
By Theorem~\ref{cor:bijection_TSPPtoCSPP} and Theorem~\ref{thm:bij_domino},
Conjecture~\ref{conj:VS} 
(Conjecture~6 of \cite{MRR2})
reduce to prove the determinantal formula in \thetag{\ref{eq:R_o}} of Conjecture~\ref{conj:result}.
We can prove all these identities in a weak form where $t=1$
by reducing the determinants to the Andrews-Burge Theorem
(Lemma~\ref{lem:AB}).
Here we use the notation 
$
(A)_{j}=A(A+1)\cdots(A+j-1).
$
%
%
%
\begin{lemma}
\label{lem:AB}
(Andrews-Burge \cite{AB})
Let
\begin{equation}
M_{n}(x,y)=\det\left(
\binom{i+j+x}{2i-j}+\binom{i+j+y}{2i-j}
\right)_{0\leq i,j\leq n-1}.
\end{equation}
Then
\begin{equation}
M_{n}(x,y)=\prod_{k=0}^{n-1}\Delta_{2k}(x+y),
\label{eq:Andrews-Burge}
\end{equation}
where $\Delta_{0}(u)=2$ and for $j>0$
\begin{equation}
\Delta_{2j}(u)=\frac{(u+2j+2)_{j}(\frac12u+2j+\frac32)_{j-1}}{(j)_{j}(\frac12u+j+\frac32)_{j-1}}.
\end{equation}
Especially,
when $y=x$,
we obtain the identity
\begin{equation}
m_{n}(x)=\det\left(
\binom{i+j+x}{2i-j}
\right)_{0\leq i,j\leq n-1}=\frac1{2^n}\prod_{k=0}^{n-1}\Delta_{2k}(2x)
\label{eq:MRR}
\end{equation}
(see \cite{A1,A2,A3,MRR3}).
\end{lemma}
\begin{theorem}
\label{thm:result}
Let $n$ be a positive integer,
and put
$r=n_{0}=\left\lceil\frac{n}{2}\right\rceil$.
\begin{enumerate}
\item[(i)]
The number of elements of $\RDPP{n}$ is equal to 
$
A_{2r+1}^\text{VS}
$
if $n=2r-1$ is odd,
i.e. $\det R_{r}^\text{o}(1)=A_{2r+1}^\text{VS}$.
\item[(ii)]
The number of elements of $\RDPP{n}$ is equal to 
$
\genfrac{}{}{}{}{(3r+2)!(2r+1)!(2r)!}{(4r+2)!(r+1)!(r!)^2}
A_{2r+1}^\text{VS}
$
if $n=2r$ is even,
i.e. 
$\det R_{r}^\text{e}(1)
=
\genfrac{}{}{}{}{(3r+2)!(2r+1)!(2r)!}{(4r+2)!(r+1)!(r!)^2}
A_{2r+1}^\text{VS}$.
\item[(iii)]
The number of elements of $\CDPP{n}$ is equal to $A_{n}^\text{HTS}$,
i.e.  $\det C_{r}^\text{o}(1)=A_{2r-1}^\text{HTS}$
and $\det C_{r}^\text{e}(1)=A_{2r}^\text{HTS}$.
\end{enumerate}
\end{theorem}
\begin{demo}{Proof}
If we put $t=1$ into $R_{r}^\text{o}(t)$ and $R_{r}^\text{e}(t)$,
then we obtain
\[
\det R_{r}^\text{o}(1)=\det_{0\leq i,j\leq r-1}\left(
\binom{i+j+1}{2i-j+1}
\right)
=\genfrac{}{}{}{}{1}{2^{r-1}}
\prod_{k=1}^{r-1}\genfrac{}{}{}{}
{\left(2k+4\right)_{k}\left(2k+\frac52\right)_{k-1}}
{\left(k\right)_{k}\left(k+\frac52\right)_{k-1}},
\]
and
\[
\det R_{r}^\text{e}(1)=\det_{0\leq i,j\leq r-1}\left(
\binom{i+j+2}{2i-j+2}
\right)
=\genfrac{}{}{}{}{1}{2^{r-1}}
\prod_{k=1}^{r-1}\genfrac{}{}{}{}
{\left(2k+6\right)_{k}\left(2k+\frac72\right)_{k-1}}
{\left(k\right)_{k}\left(k+\frac72\right)_{k-1}}
\]
from \thetag{\ref{eq:MRR}}.
A direct computation shows that
these products are equal to $A_{2r+1}^\text{VS}$
and
$
\genfrac{}{}{}{}{(3r+2)!(2r+1)!(2r)!}{(4r+2)!(r+1)!(r!)^2}
A_{2r+1}^\text{VS}
$
respectively.
Similarly,
if we put $t=1$ into $C_{r}^\text{o}(t)$ and $C_{r}^\text{e}(t)$,
then we obtain
\[
\det C_{r}^\text{o}(1)=\det_{0\leq i,j\leq r-1}\left(
\binom{i+j}{2i-j}
+\binom{i+j-1}{2i-j-1}
\right)
=\prod_{k=1}^{r-1}\genfrac{}{}{}{}
{\left(2k+1\right)_{k}\left(2k+1\right)_{k-1}}
{\left(k\right)_{k}\left(k+1\right)_{k-1}},
\]
and
\[
\det C_{r}^\text{e}(1)=\det_{0\leq i,j\leq r-1}\left(
\binom{i+j}{2i-j}
+\binom{i+j+1}{2i-j+1}
\right)
=2\prod_{k=1}^{r-1}\genfrac{}{}{}{}
{\left(2k+3\right)_{k}\left(2k+2\right)_{k-1}}
{\left(k\right)_{k}\left(k+2\right)_{k-1}},
\]
from \thetag{\ref{eq:Andrews-Burge}}.
A direct computation shows that these products
are equal to $A_{2r-1}^\text{HTS}$ and $A_{2r}^\text{HTS}$
 respectively.
\end{demo}
Thus Theorem~\ref{thm:result}(i) proves Conjecture~6 of \cite{MRR2} is true
(i.e. Conjecture~\ref{conj:VS}) when $t=1$.

{\begin{spacing}{0.5}

\end{spacing}}

\end{document}